\documentclass[aap,secthm,MSNbibl,dvips]{arximspdf}

\doi{10.1214/07-AAP513}
\volume{18}
\issue{5}
\pubyear{2008}
\firstpage{1706}
\lastpage{1736}

\makeatletter
\def\sd{\iffalse}
\newcommand{\eqref}[1]{(\ref{#1})}

\newtheorem{alem}[thm]{Lemma}
\newtheorem{aprop}[thm]{Proposition}
\newtheorem{acor}[thm]{Corollary}
\newproclaim{arem}[thm]{Remark}
\newproclaim{ex}[thm]{Example}
\newproclaim{assumption}{Assumption}

\newcommand{\E}{\mathbb{E}}

\newcommand{\N}{\mathbb{N}}
\newcommand{\dR}{\mathbb{R}}

\makeatother

\begin{document}
\begin{frontmatter}

\title{Propagation of chaos and Poincar\'e inequalities\\ for a system
of particles interacting\\ through their cdf}
\pdftitle{Propagation of chaos and Poincare inequalities\\ for a system
of particles interacting\\ through their cdf}
\runtitle{Interacting particle system through their cdf}

\begin{aug}
\author[A]{\fnms{Benjamin} \snm{Jourdain}\ead[label=e1]{jourdain@cermics.enpc.fr}}
\and
\author[B]{\fnms{Florent} \snm{Malrieu}\corref{}
\ead[label=e2]{florent.malrieu@univ-rennes1.fr}}
\runauthor{B. Jourdain and F. Malrieu}
\affiliation{\'Ecole des Ponts and Universit\'e Rennes 1}
\address[A]{CERMICS\\
\'Ecole des Ponts, ParisTech\\
6-8 av Blaise Pascal\\
Cit\'e Descartes, Champs sur Marne\\
77455 Marne-la-Vall\'ee Cedex 2\\
France\\
\printead{e1}} 
\address[B]{IRMAR\\
Universit\'e Rennes1\\
Campus de Beaulieu\\
35042 Rennes Cedex\\
France\\
\printead{e2}}
\end{aug}

\received{\smonth{1} \syear{2007}}
\revised{\smonth{12} \syear{2007}}

%
\begin{abstract}
In this paper, in the particular case of a concave flux function, we are
interested in the long time behavior of the nonlinear process
associated in
[\textit{Methodol. Comput. Appl. Probab.} \textbf{2}
(2000) 69--91] to the one-dimensional viscous scalar conservation law. We
also consider the particle system obtained by replacing the cumulative
distribution function in the drift coefficient of this nonlinear
process by
the empirical cumulative distribution function. We first obtain a
trajectorial propagation of chaos estimate which strengthens the weak
convergence result obtained in \cite{jou2000} without any convexity
assumption on the flux function. Then Poincar\'e inequalities are used to
get explicit estimates concerning the long time behavior of both the
nonlinear process and the particle system.
\end{abstract}



\pdfkeywords{65C35, 60K35, 60E15, 35K15, 46N30, Viscous scalar conservation law,
nonlinear process, particle system, propagation of chaos, Poincare
inequality, long time behavior}
\end{frontmatter}

\section*{Introduction}
In this paper, we are interested in the viscous scalar conservation law
with $C^1$ flux function $-A$
%
\begin{equation}\label{viscons}
\partial_t F_t(x)=\frac{\sigma^2}{2}\partial_{xx}F_t(x)+\partial
_x(A(F_t(x)),\qquad F_0(x)=H*m(x),
\end{equation}
where $m$ is a probability measure on the real line and $H(x)=1_{\{
x\geq0\}}$
denotes the Heaviside function. As a consequence, $H*m$ is the cumulative
distribution function of the probability measure $m$. Since $A$ appears in
this equation through its derivative, we suppose without restriction that
$A(0)=0$. According to \cite{jou2000}, one may associate the following
nonlinear process with the conservation law:
%
\begin{equation}
\cases{
X_t=X_0+\sigma B_t-\displaystyle\int_0^tA'\bigl(H*P_s(X_s)\bigr)\,ds, \cr
\forall t\geq0, \mbox{the law of $X_t$ is $P_t$},
}
\label{edsnonlin}
\end{equation}
where $(B_t)_{t\geq0}$ is a real Brownian motion independent from the
initial random variable $X_0$ with law $m$ and $\sigma$ a positive constant.
The process $X$ is said to be nonlinear in the sense that the drift
term of
the SDE depends on the entire law $P_t$ of $X_t$. More precisely, according
to \cite{jou2000}, this nonlinear stochastic differential equation
admits a
unique weak solution. Moreover, $H*P_t(x)$ is the unique bounded weak
solution of \eqref{viscons}. For $t>0$, by the Girsanov theorem, $P_t$
admits a
density $p_t$ with respect to the Lebesgue measure on the real
line.

We want to address the long time behavior of the nonlinear process
solving~\eqref{edsnonlin} by studying convergence of the density $p_t$
(see \cite{brtv} and \cite{brv} for a similar study in a different setting).
Since the cumulative distribution function $x\rightarrow H*P_s(x)$ which appears
in the
drift coefficient is nondecreasing, convexity of $A$ is a natural assumption
in order to ensure ergodicity. Then the flux function $-A$ in
the conservation law \eqref{viscons} is concave.

In the first section of the paper, after recalling results obtained in
\cite{jou2000}, we show that trajectorial uniqueness holds for
\eqref{edsnonlin} under convexity of $A$. Then we introduce a simulable
system of $n$ particles obtained by replacing in the drift coefficient the
cumulative distribution function by its empirical version and the derivative
$A'$ by a suitable finite difference approximation. When $A$ is convex,
existence and trajectorial uniqueness hold for this system. Moreover, we
prove a trajectorial estimation of propagation of chaos which
strengthens the
weak convergence result obtained in \cite{jou2000}. Unfortunately, because
the empirical cumulative distribution function is a step function and
therefore not an increasing one, this estimation is not uniform in
time.

The second and main section deals with the long time behavior of both the
nonlinear process and the particle system. We address the convergence
of the
density $p_t$ of $X_t$ by first studying the convergence of the associated
solution $H*p_t$ of \eqref{viscons} to the solution $F_\infty$ with
the same
expectation of the stationary equation
$\frac{\sigma^2}{2}\partial_{xx}F_\infty(x)+\partial_x(A(F_\infty(x))=0$
obtained by removing the time derivative in~\eqref{viscons}. For this result,
no convexity hypothesis is made on $A$. Instead, one assumes $A(u)<0$ for
$u\in(0,1)$, $A'(0)<0$, $A(1)=0$ and $A'(1)>0$. In contrast, to prove
exponential convergence of the density of the particle system uniform
in the
number $n$ of particles, we suppose that the function $A$ is uniformly
convex. This hypothesis ensures the existence of an invariant distribution
for the particle system. In \cite{soumik}, a necessary and sufficient
condition on the drift sequence is established for existence of the invariant
measure and convergence in total variation norm for the law of the particle
system at time $t$ to this measure. In the present paper, the key step to
derive quantitative convergence to equilibrium consists in obtaining a
Poincar\'e inequality for the stationary density of the particle system
uniform in $n$. This density has exponential-like tails and therefore does
not satisfy a logarithmic Sobolev inequality. So the derivation of the
Poincar\'e inequality cannot rely on the curvature criterion, used, for
instance, in \cite{cmv,cgm,mal1} or \cite{mal2} for
the granular media equation. Instead we make a direct estimation of the Poincar\'e
constant using the specific analytic form of the invariant density. To our
knowledge, our study provides the first example of a particle system, for
which a Poincar\'e inequality but no logarithmic Sobolev inequality holds
uniformly in the number $n$ of particles.

\begin{assumption*}
Throughout the paper, we assume that $A$ is a $C^1$
function on $[0,1]$ s.t. $A(0)=0$.
\end{assumption*}

\section{Propagation of chaos}
\subsection{The nonlinear process}

Let us first state existence and uniqueness for the nonlinear stochastic
differential equation \eqref{edsnonlin}.

\begin{thm}\label{thm:edsnonlin}
The nonlinear stochastic differential equation \eqref{edsnonlin}
admits a
unique weak solution $((X_t,P_t))_{t\geq0}$. For $t>0$, $P_t$ admits a
density $p_t$ with respect to the Lebesgue measure on $\dR$. The function
$(t,x)\mapsto H*P_t(x)$ is the unique bounded weak solution of the viscous
scalar conservation law \eqref{viscons}. Moreover,
%
\begin{equation} \forall
t\geq0\qquad X_t-X_0 \mbox{ is integrable and } {\mathbb E}(X_t-X_0)=-A(1)t.
\label{evolesp}
\end{equation}
Last, if the function $A$ is convex on $[0,1]$, \eqref{edsnonlin}
admits a
unique strong solution.\label{exunedsnonlin}
\end{thm}
\begin{pf}
The first and third statements are consequences of Proposition 1.2 and
Theorem 2.1 of \cite{jou2000} [uniqueness follows from uniqueness for
\eqref{viscons} and existence is obtained by a propagation of chaos
result].

According to the Yamada--Watanabe theorem, to deduce the last
statement, it is enough to check that when $A$ is convex, then trajectorial
uniqueness holds for the standard stochastic differential equation
\[
dX_t=\sigma \,dB_t-A'\bigl(H*Q_t(X_t)\bigr)\,dt
\]
where $(Q_t)_{t\geq0}$ is the flow of time-marginals of a
probability measure $Q$ on $C([0,+\infty),\dR)$. Since for each
$t\geq0$ the function $x\mapsto A'(H*Q_t(x))$ is nondecreasing, if
$(X_t)_{t\geq0}$ and $(Y_t)_{t\geq0}$ both solve this standard
SDE, then $|X_t-Y_t|$ is bounded by
\[
|X_0-Y_0|+\int_0^t{\rm
sign}(X_s-Y_s)\bigl(A'\bigl(H*Q_s(Y_s)\bigr)-A'\bigl(H*Q_s(X_s)\bigr)\bigr)\,ds,
\]
and then by $|X_0-Y_0|$ which concludes the proof of trajectorial
uniqueness.

Existence of the density $p_t$ for $t>0$ follows from the boundedness of
the drift coefficient and the Girsanov theorem. To prove \eqref{evolesp}, one
first remarks that by boundedness of the drift coefficient, for each
$t\geq
0$, the random variable $X_t-X_0$ is integrable and
\begin{eqnarray*}
{\mathbb E}(X_t-X_0)&=&-\int_0^t{\mathbb
E}\bigl(A'\bigl(H*P_s(X_s)\bigr)\bigr)\,ds\\
&=&-\int_0^t\int_{\mathbb
R}A'\bigg(\int_{-\infty}^xP_s(dy)\bigg)P_s(dx)\,ds.
\end{eqnarray*}
For $s>0$, since by the Girsanov theorem $P_s$ does not weight points,
\[
\int_{\mathbb
R}A'\biggl(\int_{-\infty}^xP_s(dy)\biggr)P_s(dx)=\bigl[A\bigl(H*P_s(x)\bigr)\bigr]_{-\infty
}^{+\infty}=A(1).
\]\upqed
\end{pf}
\begin{acor}\label{regusol}
Assume that $A$ is $C^2$ on $[0,1]$. Then the function $H*P_t(x)$ is
$C^{1,2}$ on $(0,+\infty)\times\dR$ and solves \eqref{viscons} in
the classical sense on this domain.
\end{acor}
\begin{pf}
By the Girsanov theorem, for $t_0>0$, the law $P_{t_0}$ of $X_{t_0}$
admits a density with respect to the Lebesgue measure on $\dR$.
Hence $(t,x)\mapsto H*P_{t}(x)$ is a continuous function on
$(0,+\infty)\times\dR$ with values in $[0,1]$. According to
\cite{lady}, Theorem 8.1, page 495, Remark 8.1, page 495 and Theorem 2.5,
page
18, there exists a function $u$ with values in $[0,1]$, continuous
on $[0,+\infty)\times\dR$ and $C^{1,2}$ on $(0,+\infty)\times\dR$
such that
\[
\cases{
\forall x\in\dR, &\quad $u(0,x)=H*P_{t_0}(x)$, \cr
\forall(t,x)\in(0,+\infty)\times\dR, &\quad ${\partial
_tu(t,x)=\dfrac{\sigma^2}{2}\partial_{xx}u(t,x)+\partial_x(A(u(t,x)))}$.
}
\]
By the uniqueness result for bounded weak solutions of this
viscous scalar conservation law recalled in Theorem
\ref{exunedsnonlin}, $\forall t\geq t_0$, $H*P_t(x)=u(t-t_0,x)$. The
conclusion follows since $t_0$ is arbitrary.
\end{pf}

\subsection{Study of the particle system}
For $n\in\N^*$, let ${(a_n(i))}_{1\leq i\leq n}$ be a sequence of real
numbers. In this section, we are interested in the $n$-dimensional stochastic
differential equation
%
\begin{equation}\label{systpart}
dX^{i,n}_t=\sigma \,dB^i_t-a_n\Biggl(\sum_{j=1}^n1_{\{X^{j,n}_t\leq
X^{i,n}_t\}}\Biggr)\,dt,\qquad
X^{i,n}_0=X^i_0,\ 1\leq i\leq n,
\end{equation}
where $(B^i)_{i\geq1}$ are independent standard Brownian motions
independent\vspace*{2pt} from the sequence $(X^{i}_0)_{i\geq1}$ of initial random
variables.

In the next section devoted to the approximation of the nonlinear
stochastic differential equation \eqref{edsnonlin}, we will choose
$a_n(i)$ equal to the finite difference approximation
$n(A(i/n)-A((i-1)/n))$ of $A'(\frac{i}{n})$. For this particular
choice, the
nondecreasing assumption made in the following proposition is implied
by convexity of $A$.

\begin{aprop}\label{exuneds}
Assume that the sequence $(a_n(i))_{1\leq i\leq n}$ is
nondecreasing. Then the stochastic differential equation
\eqref{systpart} has a unique strong solution. Let
$(Y^{1,n}_t,\ldots,Y^{n,n}_t)$ denote another solution starting from
$(Y^1_0,\ldots,Y^n_0)$ and driven by the same Brownian motion
$(B^1,\ldots,B^n)$. Then
%
\begin{equation}\label{decroisnorm}
a.s., \forall t\geq0\qquad
\sum_{i=1}^{n}(X^{i,n}_t-Y^{i,n}_t)^2\leq\sum
_{i=1}^{n}(X^i_0-Y^i_0)^2.
\end{equation}
In addition, if the initial conditions $(X^1_0,\ldots,X^n_0)$ and
$(Y^1_0,\ldots,Y^n_0)$ are s.t. a.s., $\forall i\in\{1,\ldots,n\}$,
$X^i_0<Y^i_0$ (resp. $X^i_0\leq Y^i_0$), then
%
\begin{equation}
\mbox{a.s., }\forall t\geq0, \forall
i\in\{1,\ldots,n\}\qquad
X^{i,n}_t<Y^{i,n}_t\mbox{ \textup{(}resp. $X^{i,n}_t\leq Y^{i,n}_t$\textup{)}}.\label{consord}
\end{equation}
\end{aprop}

Existence of a weak solution to \eqref{systpart} is a consequence of
the Girsanov theorem. Therefore, according to the  Yamada--Watanabe theorem, it
is enough to prove \eqref{decroisnorm} which implies trajectorial
uniqueness to obtain existence of a unique strong solution. To do so,
we will need the following lemma.
\begin{alem}\label{lem:permut}
Let $(a(i))_{1\leq i\leq n}$ and $(b(i))_{1\leq i\leq n}$ denote two
nondecreasing sequences of real numbers. Then for any permutation
$\tau\in{\mathcal S}_n$,
%
\begin{equation}
\label{compsom}
\sum_{i=1}^na(i)b(\tau(i))\leq\sum_{i=1}^n a(i)b(i).
\end{equation}
\end{alem}
\begin{pf}
For $n=2$, the result is an easy consequence of the inequality
\[
\bigl(a(2)-a(1)\bigr)\bigl(b(2)-b(1)\bigr)\geq0.
\]
For $n>2$, we define $\tau_1$ as $\tau$ if $\tau(1)=1$ and as
$\tau$ composed with the transposition between $1$ and
$\tau^{-1}(1)$ otherwise. This way, $\tau_1(1)=1$. In addition,
using the result for $n=2$, we get $\sum_{i=1}^na(i)b(\tau(i))\leq
\sum_{i=1}^n a(i)b(\tau_1(i))$.

For $2\leq j\leq n-1$, we define inductively $\tau_j$ as
$\tau_{j-1}$ if $\tau_{j-1}(j)=j$ and as $\tau_{j-1}$ composed with
the transposition between $j$ and $\tau_{j-1}^{-1}(j)$ otherwise.
This way, for $1\leq i\leq j$, $\tau_j(i)=i$. Again by the result
for $n=2$, one has
\[
\sum_{i=1}^na(i)b(\tau(i))\leq\sum_{i=1}^n a(i)b(\tau_1(i))\leq
\sum_{i=1}^n a(i)b(\tau_2(i))\leq\cdots\leq\sum_{i=1}^n
a(i)b(\tau_{n-1}(i)).
\]
We conclude by remarking that $\tau_{n-1}$ is the
identity.
\end{pf}

We are now ready to complete the proof of Proposition \ref{exuneds}.
\begin{pf*}{Proof of Proposition \protect\ref{exuneds}}
Let $(X^{1,n},\ldots,X^{n,n})$ and $(Y^{1,n},\ldots,Y^{n,n})$ denote
two solutions. The difference
\[
\sum_{i=1}^n (X^{i,n}_t-Y^{i,n}_t)^2-\sum_{i=1}^n (X^i_0-Y^i_0)^2
\]
is equal to
%
\begin{equation}
\label{compsol}
\hspace*{12pt} 2\int_0^t\sum_{i=1}^n(X^{i,n}_s-Y^{i,n}_s)\Biggl(a_n\Biggl(\sum
_{j=1}^n1_{\{Y^{j,n}_s\leq
Y^{i,n}_s\}}\Biggr)-a_n\Biggl(\sum_{j=1}^n1_{\{X^{j,n}_s\leq
X^{i,n}_s\}}\Biggr)\Biggr)\,ds.
\end{equation}
By the Girsanov theorem, for any $s>0$ the distributions of
$(X^{1,n}_s,\ldots,X^{n,n}_s)$ and $(Y^{1,n}_s,\ldots,Y^{n,n}_s)$ admit
densities w.r.t. the Lebesgue measure on ${\mathbb R}^n$ and therefore
$d{\mathbb P}\otimes ds$ a.e. the positions $X^{1,n}_s,\ldots,X^{n,n}_s$
(resp. $Y^{1,n}_s,\ldots,Y^{n,n}_s$) are distinct and there is a unique
permutation $\tau^X_s\in{\mathcal S}_n$ (resp. $\tau^Y_s\in{\mathcal S}_n$)
such that $X^{\tau^X_s(1),n}_s<X^{\tau^X_s(2),n}_s<\cdots<X^{\tau^X_s(n),n}_s$
(resp. $Y^{\tau^Y_s(1),n}_s<Y^{\tau^Y_s(2),n}_s<\cdots<Y^{\tau^Y_s(n),n}_s$).
Therefore $d{\mathbb P}\otimes ds$ a.e.,
\[
\sum_{i=1}^n(X^{i,n}_s-Y^{i,n}_s)\Biggl(a_n\Biggl(\sum_{j=1}^n1_{\{
Y^{j,n}_s\leq
Y^{i,n}_s\}}\Biggr)-a_n\Biggl(\sum_{j=1}^n1_{\{X^{j,n}_s\leq
X^{i,n}_s\}}\Biggr)\Biggr)
\]
is equal to
\[
\sum_{i=1}^na_n(i)\bigl(\bigl(X^{\tau^Y_s(i),n}_s-Y^{\tau
^Y_s(i),n}_s\bigr)-\bigl(X^{\tau^X_s(i),n}_s-Y^{\tau^X_s(i),n}_s\bigr)\bigr).
\]
The sequence $(a_n(i))_{1\leq i\leq n}$ is nondecreasing. Applying Lemma
\ref{lem:permut} with $b(i)=X^{\tau^X_s(i),n}_s$ and
$\tau=(\tau^X_s)^{-1}\circ\tau^Y_s$ then with $b(i)=Y^{\tau
^Y_s(i),n}_s$ and
$\tau=(\tau^Y_s)^{-1}\circ\tau^X_s$, one obtains that the integrand in
\eqref{compsol} is nonpositive $d{\mathbb P}\otimes ds$ a.e. Hence
\eqref{decroisnorm} holds.

 Let us now suppose that a.s. $\forall
i\in\{1,\ldots,n\}$, $X^i_0<Y^i_0$ and define $\nu=\inf\{t>0\dvtx\exists
i\in\{1,\ldots,n\},X^{i,n}_t\geq Y^{i,n}_t\}$ with the convention
$\inf\varnothing=+\infty$. From now on, we restrict ourselves to the event
$\{\nu<+\infty\}$. Let $i\in\{1,\ldots,n\}$ be such that
$Y^{i,n}_{\nu}=X^{i,n}_{\nu}$. There is an increasing sequence
$(s_k)_{k\geq
1}$ of\break positive times with limit $\nu$ such that $\forall k\geq
1, a_n(\sum_{j=1}^n1_{\{X^{j,n}_{s_k}\leq
X^{i,n}_{s_k}\}})<\break a_n(\sum_{j=1}^n1_{\{Y^{j,n}_{s_k}\leq
Y^{i,n}_{s_k}\}})$. Since $(a_n(i))_{1\leq i\leq n}$ is
nondecreasing, by extracting a subsequence still denoted by
$(s_k)_{k}$ for
simplicity, one deduces the existence of $j\in\{1,\ldots,n\}$ with
$j\neq i$
such that $\forall k\geq1, X^{i,n}_{s_k}<X^{j,n}_{s_k}$ and
$Y^{j,n}_{s_k}\leq Y^{i,n}_{s_k}$. Since $s_k<\nu$,
$X^{i,n}_{s_k}<X^{j,n}_{s_k}<Y^{j,n}_{s_k}\leq Y^{i,n}_{s_k}$. By continuity
of the paths, one obtains
$X^{i,n}_{\nu}=X^{j,n}_{\nu}=Y^{j,n}_{\nu}=Y^{i,n}_{\nu}$. Now since the
probability of the event
\[
\exists i_1,i_2,i_3\mbox{ dist. in
}\{1,\ldots,n\}, \exists t>0\qquad X^{i_1}_0+\sigma
B^{i_1}_t=X^{i_2}_0+\sigma B^{i_2}_t=X^{i_3}_0+\sigma
B^{i_3}_t
\]
is equal to 0, the Girsanov theorem implies that a.s. $\forall
l\in\{1,\ldots,n\}\setminus\{i,j\}$, $X^{l,n}_\nu\neq
X^{i,n}_\nu=X^{j,n}_\nu$. In the same way, $Y^{l,n}_\nu\neq
Y^{i,n}_\nu=Y^{j,n}_\nu$. By continuity of the paths and definition of
$\nu$
one deduces that for $k$ large enough, and for every $t\in[s_k,\nu]$,
\[
\mathop{\sum_{l=1}}_{l\neq
i,j}^n1_{\{Y^{l,n}_t\leq Y^{i,n}_t\}}\leq
\mathop{\sum_{l=1}}_{l\neq i,j}^n1_{\{X^{l,n}_t\leq
X^{i,n}_t\}}\mbox{; }\mathop{\sum_{l=1}}_{l\neq
i,j}^n1_{\{Y^{l,n}_t\leq Y^{j,n}_t\}}\leq
\mathop{\sum_{l=1}}_{l\neq i,j}^n1_{\{X^{l,n}_t\leq X^{j,n}_t\}}.
\]
Since a.s. $dt$ a.e., $Y^{i,n}_t\neq Y^{j,n}_t$ and
$(a_n(i))_{1\leq i\leq n}$ is nondecreasing, one obtains that a.s.
$dt$ a.e. on $[s_k,\nu]$,
\begin{eqnarray*}
&&a_n\Biggl(\sum_{l=1}^n1_{\{Y^{l,n}_t\leq
Y^{i,n}_t\}}\Biggr)+a_n\Biggl(\sum_{l=1}^n1_{\{Y^{l,n}_t\leq
Y^{j,n}_t\}}\Biggr)\\
&&\qquad\leq
a_n\Biggl(\sum_{l=1}^n1_{\{X^{l,n}_t\leq
X^{j,n}_t\}}\Biggr)+a_n\Biggl(\sum_{l=1}^n1_{\{X^{l,n}_t\leq
X^{i,n}_t\}}\Biggr).
\end{eqnarray*}
By integration with respect to $t$ on $[s_k,\nu]$, this implies
that a.s. $Y^{i,n}_\nu-X^{i,n}_\nu+Y^{j,n}_\nu-X^{j,n}_\nu\geq
Y^{i,n}_{s_k}-X^{i,n}_{s_k}+Y^{j,n}_{s_k}-X^{j,n}_{s_k}>0$.
Therefore ${\mathbb P}(\nu<+\infty)=0$.

When a.s. for $i\in\{1,\ldots,n\}$, $X^{i}_0\leq Y^i_0$, one obtains
that for $\varepsilon>0$ the solution
$(Y^{1,n,\varepsilon}_t,\ldots,Y^{n,n,\varepsilon}_t)$ to
\eqref{systpart} starting from
$(Y^1_0+\varepsilon,\ldots,Y^n_0+\varepsilon)$ is such that
\[
\mbox{a.s., }\forall t\geq0\  \forall
i\in\{1,\ldots,n\}\qquad X^{i,n}_t<Y^{i,n,\varepsilon}_t.
\]
Since by \eqref{decroisnorm}, $Y^{i,n,\varepsilon}_t\leq
Y^{i,n}_t+\sqrt{n}\varepsilon$, one easily concludes by letting
$\varepsilon\rightarrow0$.
\end{pf*}

\subsection{Trajectorial propagation of chaos}
{F}rom now on, we set
%
\begin{equation}\label{defan}
\forall n\in\N^*, \forall i\in\{1,\ldots,n\}\qquad
a_n(i)=n\biggl(A\biggl(\frac{i}{n}\biggr)-A\biggl(\frac{i-1}{n}
\biggr)\biggr)
\end{equation}
and assume that the initial positions $(X^i_0)_{i\geq1}$ of the particles
are independent and identically distributed according to $m$. We prefer
to define $a_n(i)$ with the above finite difference approximation of the
choice $A'(i/n)$ made in \cite{jou2000} because the sum $\sum_{i=1}^n
a_n(i)$ which plays a role in the long time behavior of the particle system
is then simply equal to $nA(1)$. One could also obtain trajectorial
propagation of chaos estimates similar to Theorem \ref{chaos} below for the
choice $a_n(i)=A'(i/n)$.

In the
present section, we also suppose that $A$ is a convex function on
$[0,1]$. By
Theorem \ref{thm:edsnonlin}, for each $i\geq1$, the nonlinear stochastic
differential equation
%
\begin{equation}\label{eq:nonlini}
\cases{
X^{i}_t=X^i_0+\sigma B^i_t-\displaystyle\int_0^tA'\bigl(H*P_s(X^i_s)\bigr)\,ds,\vspace*{2pt} \cr
\forall t\geq0, \mbox{the law of $X^{i}_t$ is $P_t$,}
}
\end{equation}
has a unique solution and for all $t\geq0$, the law $P_t$ of $X^i_t$
does not depend on $i$. Under a Lipschitz regularity assumption on
$A'$, we obtain the following trajectorial propagation of chaos
estimation.

\begin{thm}\label{chaos}
If $A\dvtx[0,1]\rightarrow{\mathbb R}$ is convex and $A'$ is Lipschitz
continuous with constant $K$, then
\[
\forall n\geq1, \forall1\leq i\leq n, \forall t\geq
0\qquad{\mathbb E}\biggl(\sup_{s\in
[0,t]}(X^{i,n}_s-X^i_s)^2\biggr)\leq
\frac{K^2t^2}{6n}.
\]
\end{thm}
\begin{pf}
Let us write $\sum_{i=1}^n(X^{i,n}_t-X^i_t)^2$ as
\begin{eqnarray*}
&&2\int_0^t
\sum_{i=1}^n(X^{i,n}_s-X^i_s)\Biggl(a_n\Biggl(\sum_{j=1}^n1_{\{X^j_s\leq
X^i_s\}}\Biggr)-a_n\Biggl(\sum_{j=1}^n1_{\{X^{j,n}_s\leq
X^{i,n}_s\}}\Biggr)\Biggr)\,ds\\
&&\qquad{}+2\int_0^t\sum_{i=1}^n(X^{i,n}_s-X^i_s) C(s,X^1_s,\ldots
,X^n_s)\,ds
\end{eqnarray*}
where $C(s,X^1_s,\ldots,X^n_s)$ is equal to
\[
A'\bigl(H*P_s(X^i_s)\bigr)-n\Biggl(A\Biggl(\frac{1}{n}\sum_{j=1}^n1_{\{X^j_s\leq
X^i_s\}}\Biggr)-A\Biggl(\frac{1}{n}\sum_{j=1}^n1_{\{X^j_s\leq
X^i_s\}}-\frac{1}{n}\Biggr)\Biggr).
\]
Like in the proof of trajectorial uniqueness for \eqref{systpart},
because of the convexity of $A$, the first term of the r.h.s. is
nonpositive. Moreover, by Lipschitz continuity of $A'$,
\begin{eqnarray*}
&&\Biggl(A'\bigl(H*P_s(X^i_s)\bigr)-n\Biggl(A\Biggl(\frac{1}{n}\sum_{j=1}^n1_{\{
X^j_s\leq
X^i_s\}}\Biggr)-A\Biggl(\frac{1}{n}\sum_{j=1}^n1_{\{X^j_s\leq
X^i_s\}}-\frac{1}{n}\Biggr)\Biggr)\Biggr)^2\\
&&\qquad=\Biggl(\int_0^1A'\bigl(H*P_s(X^i_s)\bigr)-A'\Biggl(\frac{1}{n}\sum_{j=1}^n1_{\{X^j_s
\leq X^i_s\}}+\frac{\theta-1}{n}\Biggr)\,d\theta\Biggr)^2\\
&&\qquad\leq\frac{K^2}{n^2}\int_0^1\Biggl(\sum_{j\neq
i}\Biggl(H*P_s(X^i_s)-1_{\{X^j_s\leq
X^i_s\}}\Biggr)+\bigl(H*P_s(X^i_s)-\theta\bigr)\Biggr)^2\, d\theta.
\end{eqnarray*}
For $s>0$, as the variables $X^i_s$ are i.i.d. with common law $P_s$
which does not weight points and $H*P_s(X^i_s)$ is uniformly
distributed on $[0,1]$,
\begin{eqnarray*}
&&\int_0^1{\mathbb E}\Biggl(\Biggl(\sum_{j\neq
i}\bigl(H*P_s(X^i_s)-1_{\{X^j_s\leq
X^i_s\}}\bigr)+\bigl(H*P_s(X^i_s)-\theta\bigr)\Biggr)^2\Biggr)\,d\theta\\
&&\qquad=\sum_{j\neq i}{\mathbb E}\bigl(\bigl(H*P_s(X^i_s)-1_{\{X^j_s\leq
X^i_s\}}\bigr)^2\bigr)+\int_0^1{\mathbb E}\bigl(\bigl(H*P_s(X^i_s)-\theta\bigr)^2\bigr)\,d\theta\\
&&\qquad=(n-1){\mathbb E}\bigl(\bigl(H*P_s(X^i_s)\bigr)\bigl(1-H*P_s(X^i_s)\bigr)\bigr)+1/6\\
&&\qquad=n/6.
\end{eqnarray*}
Using the Cauchy--Schwarz inequality, one obtains
\begin{eqnarray*}
{\mathbb
E}\Biggl(\sup_{s\in[0,t]}\sum_{i=1}^n(X^{i,n}_s-X^i_s)^2\Biggr)&\leq&
2\int_0^t\sqrt{\frac{K^2}{6n}{\mathbb
E}\Biggl(\Biggl(\sum_{i=1}^n(X^{i,n}_s-X^i_s)\Biggr)^2\Biggr)}\,ds\\
&\leq&
\frac{2K}{\sqrt{6}}\int_0^t\sqrt{{\mathbb
E}\Biggl(\sup_{u\in[0,s]}\sum_{i=1}^n(X^{i,n}_u-X^i_u)^2\Biggr)}\,ds.
\end{eqnarray*}
By comparison with the ordinary differential equation
$\alpha'(t)=2K\sqrt{\frac{\alpha(t)}{6}}$, one concludes that
\[
\forall t\geq0\qquad{\mathbb
E}\Biggl(\sup_{s\in[0,t]}\sum_{i=1}^n(X^{i,n}_s-X^i_s)^2\Biggr)\leq
\frac{K^2t^2}{6}.
\]
Exchangeability of the couples ${{((X^{i,n},X^i))}}_{i\in\{1,\ldots
,n\}}$
completes the proof.
\end{pf}

\begin{arem}\label{remalp}
One could think that assuming that $A$ is uniformly convex:
%
\begin{equation}
\exists\alpha>0, \forall0\leq x\leq y\leq1\qquad
A'(y)-A'(x)\geq\alpha(y-x)\label{unifconv}
\end{equation}
would lead
to a better estimation. Indeed, then for every $i\in\{1,\ldots,n-1\}$,
\[
a_n(i+1)-a_n(i)=n\int_{i/n}^{(i+1)/n}\biggl[A'(x)-A'\biggl(x-\frac
{1}{n}\biggr)\biggr]\,dx\geq
\frac{\alpha}{n}.
\]
But since even in this situation, the nonpositive term
\[
\sum_{i=1}^n(X^{i,n}_s-X^i_s)\Biggl(a_n\Biggl(\sum_{j=1}^n1_{\{X^j_s\leq
X^i_s\}}\Biggr)-a_n\Biggl(\sum_{j=1}^n1_{\{X^{j,n}_s\leq
X^{i,n}_s\}}\Biggr)\Biggr)
\]
vanishes as soon as the order between the coordinates of
$(X^{1,n}_s,\ldots,X^{n,n}_s)$ is the same as the order between the
coordinates of $(X^{1}_s,\ldots,X^{n}_s)$, we were not able so far
to improve the estimation.
\end{arem}
\begin{acor}\label{cor:nonlin} Under the hypotheses of Theorem \ref
{chaos}, let $\tilde{m}$ be a probability measure on $\mathbb R$ such
that $\forall x\in{\mathbb R}$, $H*\tilde{m}(x)\leq H*m(x)$. If for
some random variable $U_1$ uniform on $[0,1]$ independent from
$(B^i)_{i\geq1}$, $X^1_0=\inf\{x\dvtx H*m(x)\geq U_1\}$ and
$(Y^1_t)_{t\geq0}$ denotes the solution of the nonlinear stochastic
differential equation
%
\begin{equation}
\cases{Y^{1}_t=Y^1_0+\sigma B^1_t-\displaystyle\int_0^tA'\bigl(H*\tilde
{P}_s(Y^1_s)\bigr)\,ds, \cr
\forall t\geq0, \mbox{the law of $Y^1_t$ is
$\tilde{P}_t$},
}
\end{equation}
with $Y^1_0=\inf\{x\dvtx H*\tilde{m}(x)\geq U_1\}$, then
\[
{\mathbb
P}(\forall t\geq0, X^1_t\leq Y^1_t)=1.
\]
Moreover $\forall t\geq
0$, $\forall x\in{\mathbb R}$, $H*\tilde{P}_t(x)\leq H*P_t(x)$.
Last, the function $t\mapsto{\mathbb E}|Y^1_t-X^1_t|$ is constant.
\end{acor}
\begin{arem}
At least when $m$ and $\tilde{m}$ do not weight points, one has a.s.
$A'(H*P_0(X^1_0))=A'(H*\tilde{P}_0(Y^1_0))$ since
$H*m(X^1_0)=H*\tilde{m}(Y^1_0)=U_1$. Therefore a.s. $d(Y^1-X^1)_0=0$
and one may wonder whether a.s. $Y^1_t-X^1_t$ does not depend on
$t$. If this property holds, necessarily, a.s. $dt$ a.e.
$A'(H*P_t(X^1_t))=A'(H*\tilde{P}_t(Y^1_t))$. If $A'$ is increasing,
a.s. for all $t>0$, $H*p_t(X^1_t)=H*\tilde{p}_t(Y^1_t)$ with $p_t$
and $\tilde{p}_t$ denoting the respective densities of $P_t$ and
$\tilde{P}_t$. If $A$ is $C^2$, the Brownian contribution in
$d(H*p_t(X^1_t)-H*\tilde{p}_t(Y^1_t))$ given by It\^o's
formula vanishes, that is, $p_t(X^1_t)=\tilde{p}_t(Y^1_t)$ and $\forall
u\in\,]0,1[$,
$p_t((H*p_t)^{-1}(u))=\tilde{p}_t((H*\tilde{p}_t)^{-1}(u))$ or
equivalently $((H*p_t)^{-1})'(u)=((H*\tilde{p}_t)^{-1})'(u)$. Hence
$Y^1_t=X^1_t+c$ for a deterministic constant $c$ which does not
depend on $t$ according to \eqref{evolesp}. Letting $t\rightarrow
0$, one obtains $Y^1_0=X^1_0+c$. This necessary condition turns out
to be sufficient as $(X^1_t+c)_{t\geq0}$ obviously solves the
nonlinear stochastic differential equation \eqref{edsnonlin}
starting from $X^1_0+c$.
\end{arem}
\begin{pf*}{Proof of Corollary \protect\ref{cor:nonlin}}
For $(U_i)_{i\geq2}$ a sequence of\vspace*{2pt} independent uniform random
variables independent from $(U_1,(B^i)_{i\geq1})$, we set
\[
\forall i\geq2\qquad X^i_0=\inf\{x\dvtx H*m(x)\geq U_i\}\quad\mbox{and}\quad
Y^i_0=\inf\{x\dvtx H*\tilde{m}(x)\geq U_i\}.
\]
Since $H*\tilde{m}\leq H*m$, a.s. $\forall i\geq1$, $Y^i_0\geq
X^i_0$. From Proposition \ref{exuneds}, one deduces that the
solutions $(X^{1,n}_t,\ldots,X^{n,n}_t)$ and
$(Y^{1,n}_t,\ldots,Y^{n,n}_t)$ to \eqref{systpart} respectively
starting from $(X^{1}_0,\ldots,X^{n}_0)$ and
$(Y^{1}_0,\ldots,Y^{n}_0)$ are such that
\[
\mbox{a.s., }\forall n\geq1, \forall i\in\{1,\ldots,n\}, \forall
t\geq0\qquad Y^{i,n}_t\geq X^{i,n}_t.
\]
Since, by Theorem \ref{chaos}, for fixed $t\geq0$, one may
extract from $(X^{1,n}_t,Y^{1,n}_t)_{n\geq1}$ a subsequence almost
surely converging to $(X^1_t,Y^1_t)$, one easily deduces that
${\mathbb P}(\forall t\geq0, X^1_t\leq Y^1_t)=1$. Hence
\[
\forall t\geq0, \forall x\in{\mathbb
R} \qquad H*\tilde{P}_t(x)={\mathbb P}(Y^1_t\leq x)\leq{\mathbb
P}(X^1_t\leq x)=H*P_t(x).
\]
Since $|Y^1_t-X^1_t|-|Y^1_0-X^1_0|=Y^1_t-Y^1_0-(X^1_t-X^1_0)$, \eqref
{evolesp} ensures that ${\mathbb
E}|Y^1_t-X^1_t|\in[0,+\infty]$ does not depend on $t$.
\end{pf*}

\section{Long time behavior}\label{sec:ltb}
In this section we are interested in the long time behavior of both the
nonlinear process and the particle system. According to \eqref{evolesp} and
the equality $\sum_{i=1}^n a_n(i)=nA(1)$ which follows from \eqref
{defan}, we have to suppose $A(1)=0$ in order
to obtain convergence of the densities as $t$ tends to infinity. We address
the convergence of the density $p_t$ of $X_t$ by first studying the
convergence of the associated cumulative distribution function $F_t$
under the following hypothesis denoted by \eqref{hyp} in the
sequel:
\renewcommand{\theequation}{H}
\begin{eqnarray}\label{hyp}
 A(0)&=&A(1)=0,\qquad A'(0)<0,\nonumber\\[-8pt]\\[-8pt]
 A'(1)&>&0 \quad\mbox{and}\quad\forall u\in(0,1)\qquad A(u)<0.\nonumber
\end{eqnarray}
These assumptions determine the spatial behavior at infinity of the drift
coefficient in \eqref{edsnonlin}.

\renewcommand{\theequation}{\arabic{equation}}
\setcounter{equation}{12}

To prove exponential convergence of the
density of the particle system uniform in the number $n$ of particles,
we make
the stronger assumption of uniform convexity on $A$. The key step in
the proof
is to obtain a Poincar\'e inequality uniform in $n$ for the stationary density
of the particle system. This density has exponential-like tails and therefore
does not satisfy a logarithmic Sobolev inequality. So the derivation of the
Poincar\'e inequality cannot rely on the curvature criterion, used, for
instance, by Malrieu \cite{mal1,mal2} when dealing with the granular
media equation. Instead, we take advantage of the following nice
feature: up
to reordering of the coordinates, the stationary density is the density
of the
image by a linear transformation of a vector of independent exponential
variables. And it turns out that the control of the constant in the
$n$-dimensional Poincar\'e inequality relies on the Hardy inequality
stated in
Lemma \ref{lemhard} which is a one-dimensional Poincar\'e-like
inequality. To
our knowledge, our study provides the first example of a particle
system, for
which a Poincar\'e inequality but no logarithmic Sobolev inequality holds
uniformly in the number $n$ of particles.

\subsection{The nonlinear process}
In this section, we are first going to obtain necessary and sufficient
conditions on the function $A$ ensuring existence for the stationary
Fokker--Planck equation obtained by removing the time-derivative in the
nonlinear Fokker--Planck equation
%
\begin{equation}\label{fp}
\partial_t
p_t=\frac{\sigma^2}{2}\partial_{xx}p_t+\partial_x\bigl(A'(H*p_t)p_t\bigr)
\end{equation}
satisfied by the density of the solution of \eqref{edsnonlin}. Under a
slightly stronger condition, the solutions satisfy a Poincar\'e
inequality.

\begin{alem}\label{fpstat}
A necessary and
sufficient condition for the existence of a probability measure $\mu$
solving the stationary Fokker--Planck equation
\[
\frac{\sigma^2}{2}\partial_{xx}\mu+\partial_x\bigl(A'\bigl(H*\mu(x)\bigr)\mu\bigr)=0
\]
in the distribution sense is $A(1)=0$ and $A(u)<0$ for all $u\in(0,1)$.
Under that condition, all the solutions are the translations of a
probability measure with a $C^1$ density $f$ which satisfies
%
\begin{eqnarray}\label{der1}
\forall x\in\dR\qquad
f(x)&=&-\frac{2}{\sigma^2}A\bigl(H*f(x)\bigr)\quad\mbox{and}\nonumber\\[-8pt]\\[-8pt]
f'(x)&=&-\frac{2}{\sigma^2}A'\bigl(H*f(x)\bigr)f(x).\nonumber
\end{eqnarray}
If $A'(0)<0$ and $A'(1)>0$, then
\begin{eqnarray}\label{equivinf}
f(x)&\sim& \cases{%
-\displaystyle\frac{2A'(0)}{\sigma^2}\int_{-\infty}^xf(y)\,dy, &\quad when $x\rightarrow-\infty$,\cr
\displaystyle\frac{2A'(1)}{\sigma^2}\int_x^{+\infty}f(y)\,dy, & \quad when $x\rightarrow+\infty$,}
\nonumber\\[-8pt]
\\[-8pt]
\int_0^x\frac{dy}{f(y)}&\sim&
\cases{\displaystyle\frac{-\sigma^2}{2A'(0)f(x)}, &\quad when $x\rightarrow-\infty$,\cr
\displaystyle\frac{\sigma^2}{2A'(1)f(x)}, & \quad when $x\rightarrow+\infty$,}
\nonumber
\end{eqnarray}
and all the solutions satisfy a Poincar\'e inequality and have a
finite expectation. Last, if the function $A$ is $C^2$ on $[0,1]$,
then $f$ is $C^2$ and satisfies
%
\begin{equation}
f''(x)=-\frac{2}{\sigma^2}A''\bigl(H*f(x)\bigr)f^2(x)+\frac{f'^2(x)}{f(x)}.\label{der2}
\end{equation}
\end{alem}
\begin{pf}
Let $\mu$ be a probability measure on ${\mathbb R}$ solving the
stationary Fokker--Planck equation. The equality
$\frac{\sigma^2}{2}\partial_{xx}\mu=-\partial_x(A'(H*\mu(x))\mu)$
ensures that $\mu$ does not weight points. Hence the stationary
equation is equivalent to
$\partial_{xx}(\frac{\sigma^2}{2}\mu+A(H*\mu(x)))=0$. One
deduces that $\mu$ possesses a $C^1$ density $f$ such that
%
\begin{equation}
\forall x\in{\mathbb R}\qquad f(x)=-\frac{2}{\sigma^2}A\bigl(H*f(x)\bigr)+\alpha
x+\beta,\label{edointeg}
\end{equation}
for some constants $\alpha$ and $\beta$. Since $A(0)=0$, letting
$x\rightarrow
-\infty$ then $x\rightarrow+\infty$ in the last equality, one obtains
$\alpha=\beta=A(1)=0$. For $u\in(0,1)$, since $u=H*f(x)$ for some $x\in
\dR$
and $H*f$ is not constant and equal to $u$, the Cauchy--Lipschitz
theorem and
\eqref{edointeg} imply that $A(u)\neq0$. Since $f$ is nonnegative, $A(u)<0$.
Hence $A(1)=0$ and $A(u)<0$ for all $u\in(0,1)$ is a necessary
condition.\par
Under that condition, a probability measure $\mu$ solves the stationary
Fokker--Planck equation if and only if its cumulative distribution function
$H*\mu(x)$ is a $C^2$ solution to the differential equation
%
\begin{equation}\label
{edofoncrep}
\varphi'(x)=-\frac{2}{\sigma^2}A(\varphi(x)),\qquad x\in{\mathbb R}.
\end{equation}
By the Cauchy--Lipschitz theorem, for each $v\in[0,1]$ this equation
admits a
unique solution $\varphi_v$ defined on $\dR$ with values in $[0,1]$
such that
$\varphi_v(0)=v$. Moreover, as $A(0)=A(1)=0$, $\varphi_0\equiv0$ and
$\varphi_1\equiv1$ and
%
\begin{equation}\label{eq:var}
\forall v\in(0,1), \forall x\in\dR\qquad 0<\varphi_v(x)<1.
\end{equation}
For $v\in(0,1)$, since $\varphi_v$ is nondecreasing and
$\varphi_v(x)=v-\frac{2}{\sigma^2}\int_0^xA(\varphi_v(y))\,dy$,
necessarily $\lim_{y\rightarrow+\infty} \varphi_v(y)=1$. In the same
way, $\lim_{y\rightarrow-\infty} \varphi_v(y)=0$ and $\varphi_v$ is
an increasing function from $\dR$ to $(0,1)$ with inverse denoted by
$\varphi_v^{-1}$. The uniqueness result for \eqref{edofoncrep} implies
that $\forall v\in(0,1), \forall x\in{\mathbb R},
\varphi_v(x)=\varphi_{{1}/{2}}(x+\varphi_{{1}/{2}}^{-1}(v))$.
Therefore the solutions to the stationary Fokker--Planck equation are
the probability measures obtained by spatial translation of the
probability measure with density $f(x)=\varphi'_{{1}/{2}}(x)$
which satisfies \eqref{der1} according to \eqref{edofoncrep}.

Let us now suppose that $A'(0)<0$ and $A'(1)>0$. When $x\rightarrow
+\infty$,
\[
f(x)=-\frac{2}{\sigma^2}A\biggl(1-\int_x^{+\infty}f(y)\,dy\biggr)\sim\frac
{2A'(1)}{\sigma^2}\int_x^{+\infty}f(y)\,dy.
\]
By \eqref{der1}, $\frac{f'(x)}{f(x)}=(\log
f(x))'=-\frac{2}{\sigma^2}A'(\varphi_{{1}/{2}}(x))$ converges to
$-\frac{2A'(1)}{\sigma^2}$ as $x\rightarrow+\infty$. This implies that
$\frac{\log(f(x))}{x}$ converges to $-\frac{2A'(1)}{\sigma^2}$ and that
$xf(x)1_{\{x\geq0\}}$ is integrable. Moreover, since
$\int_0^{+\infty}\frac{dy}{f(y)}=+\infty$,
$\int_0^x\frac{dy}{f(y)}\sim\frac{\sigma^2}{2A'(1)}\int_0^x-\frac
{f'(y)}{f^2(y)}\,dy\sim\frac{\sigma^2}{2A'(1)f(x)}$,
as $x\rightarrow
+\infty$.
In the same way, one obtains the equivalents
given in \eqref{equivinf} when $x\rightarrow-\infty$ and checks the
integrability of the function
$xf(x)1_{\{x\leq0\}}$. From \eqref{equivinf}, one has
\[
\lim_{x\rightarrow-\infty}\int_{-\infty}^x
f(y)\,dy\int_x^0\frac{dy}{f(y)}=\frac{\sigma^4}{4(A'(0))^2}
\]
and
\[
\lim_{x\rightarrow+\infty}\int_x^{+\infty}
f(y)\,dy\int_0^x\frac{dy}{f(y)}=\frac{\sigma^4}{4(A'(1))^2}.
\]
By Theorem 6.2.2, page 99 of \cite{abcfgmrs}, one concludes that the measure
with density $f$ satisfies a Poincar\'e inequality.

By \eqref{der1}, the function $f$ is $C^2$ as soon as the function $A$
is $C^2$ on $[0,1]$. Moreover,
$f''(x)=-\frac{2}{\sigma^2}A''(H*f(x))f^2(x)-\frac{2}{\sigma^2}A'(H*f(x))f'(x)$
which combined with \eqref{der1} implies \eqref{der2}.
\end{pf}
\begin{arem}
When $A$ is a $C^1$ convex function on $[0,1]$ such that
$A(0)=A(1)=0$ and $A'(u)<0$ for some $u\in(0,1)$, then the
necessary and sufficient condition in Lemma \ref{fpstat} is
obviously satisfied. Since \eqref{der1} implies
\begin{eqnarray*}
(\log
f(x))''&=&\biggl(\frac{f'(x)}{f(x)}\biggr)'=\biggl(\frac{-{2}/{\sigma
^2}A'(H*f(x))f(x)}{f(x)}\biggr)'\\
&=&-\frac{2}{\sigma^2}A''\bigl(H*f(x)\bigr)f(x)\leq
0,
\end{eqnarray*}
the probability measures solving the stationary Fokker--Planck
equation admit log-concave densities with respect to the Lebesgue
measure. Log-concavity is a property stronger than the existence of a
Poincar\'e inequality (see \cite{fou2005}).
\end{arem}
\begin{ex}
Using \eqref{edofoncrep} and \eqref{eq:var}, the following two choices for
$A$ lead to exact computations and different tails for the stationary
densities:
\begin{itemize}
\item if $A(x)=\frac{1}{2}x(x-1)$, one gets
$\log(\frac{\varphi_{{1}/{2}}(x)}{1-\varphi_{
{1}/{2}}(x)})=x/\sigma^2$, that is,
\[
\varphi_{{1}/{2}}(x)=\frac{e^{x/\sigma^2}}{1+e^{x/\sigma^2}}
\quad\mbox{and}\quad
\varphi'_{{1}/{2}}(x)=\frac{1}{4\sigma^2\cosh^2(x/2\sigma^2)};
\]
\item if $A(x)=x^3-x=x(x-1)(x+1)$,
\[
\varphi_{\sqrt{1/2}}(x)=\frac{1}{\sqrt{1+e^{-4x/\sigma^2}}}
\quad\mbox{and}\quad
\varphi'_{\sqrt{1/2}}(x)=\frac{2e^{-4x/\sigma^2}}{\sigma
^2(1+e^{-4x/\sigma^2})^{3/2}}.
\]
\end{itemize}
\end{ex}

When $A(1)=0$ and $A(u)<0$ for all $u\in(0,1)$, a natural question is
how to
link the translation parameter of the candidate long time limit of the
marginal $P_t$ solving the stationary Fokker--Planck equation to the initial
marginal $m$. When $\int_{\mathbb R}|x|m(dx)<+\infty$, by \eqref
{evolesp}, for
all $t\geq0$, ${\mathbb E}(X^1_t)={\mathbb E}(X^1_0)$. Therefore the
translation parameter is chosen in order to ensure that the invariant measure
has the same mean as the initial measure
$m$.

Let us denote by $p_t$ the density of $P_t$ and by $F_t=H*P_t$ its cumulative
distribution function.
\begin{thm}\label{th:nonlin}
Let $A$ be $C^2$ on $[0,1]$ satisfying \eqref{hyp}. Assume that $m$
admits a
density $p_0$ such that $\int_\dR|x|p_0(x)\,dx<+\infty$ and
$\int_\dR\frac{(p_0(x)-p_\infty(x))^2}{p_\infty(x)}\,dx$ is small enough where
$p_\infty$ denotes the stationary distribution with same expectation as
$p_0$. Last, we suppose that $A$ and $p_0$ are such that $p$ is a smooth
solution of \eqref{fp}. Then
$\int_\dR\frac{(p_t(x)-p_\infty(x))^2}{p_\infty(x)}\,dx$ converges to $0$
exponentially fast as $t\rightarrow+\infty$.
\label{longtime}
\end{thm}

By a smooth solution of \eqref{fp}, we mean that $p$ possesses enough
regularity and integrability so that the formal computations made in the
proof below are justified.

\begin{ex}\label{ex:explicit}
When $A(x)=\frac{1}{2}(x^2-x)$, one easily checks that the function
$\phi(t,x)=-F_t(x+\frac{t}{2})$ solves Burgers' equation
\[
\partial_t\phi=\frac{\sigma^2}{2}\partial_{xx}\phi-\frac{1}{2}\partial
_x\phi^2,\qquad \phi(0,x)=-F_0(x).
\]
By the Cole--Hopf transformation,
$\psi(t,x)=\exp(-\frac{1}{\sigma^2}\int_{-\infty}^x\phi
(t,y)\,dy)$
solves the heat equation
\[
\partial_t\psi=\frac{\sigma^2}{2}\partial_{xx}\psi,\qquad \psi(0,x)=\exp
\biggl(\frac{1}{\sigma^2}\int_{-\infty}^xF_0(y)\,dy\biggr).
\]
Since
$F_t(x)=\sigma^2\frac{\partial_x\psi}{\psi}(t,x-\frac{t}{2})$, one
deduces that
%
\begin{equation}
F_t(x)=\frac{\int_\dR e^{-{(x-{t}/{2}-y)^2}/{2\sigma^2
t}}F_0(y)\psi(0,y)\,{dy}/{(\sigma\sqrt{2\pi t})}}{\int_\dR
e^{-{(x-{t}/{2}-y)^2}/{2\sigma^2
t}}\psi(0,y)\,{dy}/{(\sigma\sqrt{2\pi t})}}.\label{rapport}
\end{equation}
If $\bar{x}$ denotes the expectation associated with the cumulative
distribution function $F_0$, one has
$\int_{-\infty}^{\bar{x}}F_0(z)\,dz=\int_{\bar{x}}^{+\infty}(1-F_0(z))\,dz$.
Since
\[
\int_{-\infty}^{x}F_0(z)\,dz=\int_{-\infty}^{\bar{x}}F_0(z)\,dz-\int_{\bar
{x}}^{x}\bigl(1-F_0(z)\bigr)\,dz+(x-\bar{x}),
\]
one deduces that the function
$\tilde{\psi}(0,x)=e^{-{(x-\bar{x})}/{\sigma^2}}\psi(0,x)$ [resp.
$\psi(0,x)$] is bounded on $\dR_+$ (resp. $\dR_-$) and converges to
$1$ as $x$ tends to $+\infty$ (resp. $-\infty$).

Let us deduce the limit of $F_t(x)$ as $t\rightarrow+\infty$. Writing
the integral for $y\in\dR$ as the sum of the integrals for $y\in
\dR_-$ and for $y\in\dR_+$, and making the change of variables
$z=\frac{y-x+{t}/{2}}{\sigma\sqrt{t}}$ (resp.
$z=\frac{y-x-{t}/{2}}{\sigma\sqrt{t}}$) in the first (resp.
second) integral, one obtains
\begin{eqnarray*}
&&\int_\dR e^{-{(y-x+{t}/{2})^2}/{(2\sigma^2
t)}}F_0(y)\psi(0,y)\frac{dy}{\sigma\sqrt{2\pi t}}\\
&&\qquad=\int_\dR
e^{-{z^2}/{2}}1_{\{z\leq
{\sqrt{t}}/{(2\sigma)}-{x}/{(\sigma\sqrt{t})}\}}\\
&&\phantom{\qquad=\int_\dR}{}\times F_0\biggl(\sigma\sqrt
{t}z+x-\frac{t}{2}\biggr)\psi\biggl(0,\sigma\sqrt{t}z+x-\frac{t}{2}\biggr)\frac{dz}{\sqrt
{2\pi}}\\
&&\qquad\quad{}+e^{{(x-\bar{x})}/{\sigma^2}}\int_\dR
e^{-{z^2}/{2}}1_{\{z\geq-{\sqrt{t}}/{(2\sigma)}-{x}/{(\sigma
\sqrt{t})}\}}\\
&&\phantom{\qquad\quad{}+e^{{(x-\bar{x})}/{\sigma^2}}\int_\dR}{}\times F_0\biggl(\sigma\sqrt{t}z+x+\frac{t}{2}\biggr)\tilde{\psi}\biggl(0,\sigma\sqrt
{t}z+x+\frac{t}{2}\biggr)\frac{dz}{\sqrt{2\pi}}.
\end{eqnarray*}
By the Lebesgue theorem, the first term of the right-hand side converges to
$0$ whereas the second term converges to
$e^{{(x-\bar{x})}/{\sigma^2}}$. Replacing $F_0$ by $1$ in the above
computation, one obtains that the denominator in \eqref{rapport}
converges to $1+e^{{(x-\bar{x})}/{\sigma^2}}$. Therefore
\[
\forall x\in\dR\qquad \lim_{t\rightarrow+\infty}
F_t(x)=\frac{e^{{(x-\bar{x})}/{\sigma^2}}}{1+e^{{(x-\bar
{x})}/{\sigma^2}}}.
\]
Notice that in the same way, one may also obtain the limit of the
density
\begin{eqnarray*}
p_t(x)&=&\frac{\int_\dR
({(y+{t}/{2}-x)}/{(\sigma^2t)})e^{-{(x-{t}/{2}-y)^2}/{(2\sigma^2
t)}}F_0(y)\psi(0,y)\,{dy}/{(\sigma\sqrt{2\pi t})}}{\int_\dR
e^{-{(x-\frac{t}{2}-y)^2}/{(2\sigma^2
t)}}\psi(0,y)\,{dy}/{(\sigma\sqrt{2\pi
t})}}\\
&&{}-\frac{1}{\sigma^2}\biggl(\frac{\int_\dR
e^{-{(x-{t}/{2}-y)^2}/{(2\sigma^2
t)}}F_0(y)\psi(0,y)\,{dy}/{(\sigma\sqrt{2\pi t})}}{\int_\dR
e^{-{(x-{t}/{2}-y)^2}/{(2\sigma^2
t)}}\psi(0,y)\,{dy}{(\sigma\sqrt{2\pi t})}}\biggr)^2.
\end{eqnarray*}
One easily checks
\begin{eqnarray*}
\forall x\in\dR\qquad \lim_{t\rightarrow
+\infty}p_t(x)&=&\frac{1}{\sigma^2}\biggl(\frac{e^{{(x-\bar{x})}/{\sigma
^2}}}{1+e^{{(x-\bar{x})}/{\sigma^2}}}-\frac{e^{{2(x-\bar
{x})}/{\sigma^2}}}{(1+e^{{(x-\bar{x})}/{\sigma^2}})^2}
\biggr)\\
&=&\frac{1}{4\sigma^2\cosh^2({(x-\bar{x})}/{2\sigma^2})}.
\end{eqnarray*}
\end{ex}

In order to prove Theorem \ref{longtime}, we are first going to check
exponential convergence of $F_t$ to the cumulative distribution function
$F_\infty$ of $p_\infty$. Let $G_t=F_t-F_\infty$. Since for a random variable
$X$ with cumulative distribution function $F$, ${\mathbb
E}(X)=\int_0^{+\infty}(1-F(x))\,dx-\int_{-\infty}^0F(x)\,dx$, the equality
of the
expectations associated to $F_t$ and $F_\infty$ writes $\int_\dR
G_t(x)\,dx=0$. This very convenient expression of the link between
$p_t$ and
$p_\infty$ is one main reason for first considering the convergence of $G_t$
to $0$. In order to prove this convergence, we need the following result.

\begin{alem}\label{lem:G}
Under the assumptions of Theorem \textup{\ref{th:nonlin}}, one has
%
\begin{equation}
\int_\dR\frac{G_t^2(x)}{p_\infty(x)}\,dx\leq c\int_\dR
{\biggl(\frac{G_t(x)}{p_\infty(x)}\biggr)'}^2p_\infty(x)\,dx \label{poincarérep}
\end{equation}
where $c$ denotes the constant in the Poincar\'e
inequality satisfied by $p_\infty$. Moreover
\begin{eqnarray}\label{egal}
&&\int_\dR\frac{(p_t(x)-p_\infty(x))^2}{p_\infty(x)}\,dx\nonumber\\[-8pt]\\[-8pt]
&&\qquad=\int_\dR
{\biggl(\frac{G_t(x)}{p_\infty(x)}\biggr)'}^2p_\infty(x)\,dx+\frac
{2}{\sigma^2}\int_\dR G_t(x)^2A''(F_\infty)(x)\,dx\nonumber
\end{eqnarray}
and
\begin{eqnarray}\label{majofoncrepdens}
\int_\dR\frac{G_t(x)^2}{p_\infty(x)}\,dx\leq
\tilde{c}\int_\dR\frac{(p_t(x)-p_\infty(x))^2}{p_\infty(x)}\,dx.
\end{eqnarray}
\end{alem}
\begin{arem}
When $A$ is convex, \eqref{majofoncrepdens} is a consequence of
\eqref{egal} and \eqref{poincarérep}.
\end{arem}

\begin{pf*}{Proof of Lemma \protect\ref{lem:G}}
As $\int_\dR G_t(x)\,dx=0$, \eqref{poincarérep} is the Poincar\'e
inequality satisfied by $p_\infty$
written for the function $G_t/p_\infty$.

Since
$(\frac{G_t(x)}{p_\infty(x)})'=\frac{G_t'(x)}{p_\infty
(x)}-\frac{G_t(x)p_\infty'(x)}{p_\infty(x)^2}$,
one has
\begin{eqnarray*}
\int_\dR
{\biggl(\frac{G_t(x)}{p_\infty(x)}\biggr)'}^2p_\infty(x)\,dx&=&\int_\dR
\frac{(p_t(x)-p_\infty(x))^2}{p_\infty(x)}\,dx-\int_\dR
\frac{{G_t^2(x)}'p_\infty'(x)}{p_\infty^2(x)}\,dx\\
&&{}+\int_\dR\frac{G_t^2(x){p_\infty'(x)}^2}{p_\infty^3(x)}\,dx\\
&=&\int_\dR\frac{(p_t(x)-p_\infty(x))^2}{p_\infty(x)}\,dx+\int_\dR
\frac{{G_t^2(x)}p_\infty''(x)}{p_\infty^2(x)}\,dx\\
&&{}-\int_\dR\frac{G_t^2(x){p_\infty'(x)}^2}{p_\infty^3(x)}\,dx.
\end{eqnarray*}
Since $p_\infty$ solves \eqref{der2}, one easily deduces \eqref{egal}.

Writing $G^2_t(y)$ as
\[
2\biggl(1_{\{y\leq0\}}\int_{-\infty}^y
G_t(p_t-p_\infty)(x)\,dx-1_{\{y>0\}}\int_y^{+\infty}
G_t(p_t-p_\infty)(x)\,dx\biggr),
\]
one obtains
\begin{equation}\label{integpart}
\int_\dR\frac{G_t^2}{p_\infty}(x)\,dx=-2\int_\dR G_t(p_t-p_\infty)(x)\int
_0^x\frac{1}{p_\infty(y)}\,dy \,dx.
\end{equation}
By \eqref{equivinf}, and since $\frac{1}{p_\infty}$ is bounded from
below and above on each
compact subset of the real line,
\[
\exists C>0,\forall
x\in\dR\qquad \bigg|\int_0^{x}\frac{1}{p_\infty(y)}\,dy\bigg|\leq
\frac{C}{p_\infty(x)}.
\]
Using the Cauchy--Schwarz inequality in \eqref{integpart}, and inserting
the latter bound, one obtains
\[
\int_\dR\frac{G_t^2}{p_\infty}(x)\,dx\leq2C\biggl(\int_\dR
\frac{G_t^2}{p_\infty}(x)\,dx\biggr)^{1/2}\biggl(\int_\dR
\frac{(p_t(x)-p_\infty(x))^2}{p_\infty(x)}\,dx\biggr)^{1/2}.
\]
One easily deduces \eqref{majofoncrepdens}.
\end{pf*}

According to \eqref{majofoncrepdens}, the exponential convergence of
$\int_\dR
\frac{(p_t(x)-p_\infty(x))^2}{p_\infty(x)}\,dx$ to zero is a stronger
result than
the exponential convergence stated in the next lemma.

\begin{alem}\label{convfoncrep}
Under the assumptions of Theorem \ref{th:nonlin}, there is a positive
constant $C$ such that if $\int_\dR\frac{G_0^2}{p_\infty}(x)\,dx$ is small
enough, then
\[
\forall t\geq0\qquad\int_\dR\frac{G_t^2}{p_\infty}(x)\,dx\leq
\frac{e^{-C t}}{C}\int_\dR\frac{G_0^2}{p_\infty}(x)\,dx.
\]
\end{alem}

\begin{pf}
According to \eqref{der1}, one has
$\frac{\sigma^2}{2}F_\infty''+(A(F_\infty))'=0$ which also writes
$\frac{p_\infty'}{p_\infty}=-\frac{2}{\sigma^2}A'(F_\infty)$.
Combining these equations with \eqref{viscons}, then using Young's
inequality, one easily obtains for $\varepsilon>0$,
\begin{eqnarray}\label{majevol}
&&\frac{1}{2}\frac{d}{dt}\int_\dR
\frac{G_t^2}{p_\infty}(x)\,dx\nonumber\\
&&\qquad=-\frac{\sigma^2}{2}\int_\dR{\biggl(\frac
{G_t(x)}{p_\infty(x)}\biggr)'}^2p_\infty(x)\,dx\nonumber\\[-8pt]\\[-8pt]
&&\qquad\quad{}-\int_\dR\bigl(A(F_t)-A(F_\infty)-A'(F_\infty)G_t\bigr)(x)\biggl(\frac
{G_t(x)}{p_\infty(x)}\biggr)'\,dx\nonumber\\
&&\qquad\leq\biggl(\varepsilon-\frac{\sigma^2}{2}\biggr)\int_\dR{\biggl(\frac{G_t}{p_\infty
}(x)\biggr)'}^2
p_\infty(x)\,dx+\frac{\|A''\|^2_\infty}{16\varepsilon}\int_\dR\frac
{G_t^4(x)}{p_\infty(x)}\,dx.\nonumber
\end{eqnarray}
Since
%
\begin{eqnarray}\label{majGinf}
\|G_t\|^2_\infty&\leq&\biggl(\int_\dR\frac{|p_t(x)-p_\infty(x)|}{\sqrt
{p_\infty(x)}}\sqrt{p_\infty(x)}\,dx\biggr)^2\nonumber\\[-8pt]\\[-8pt]
&\leq&\int_\dR
\frac{(p_t(x)-p_\infty(x))^2}{p_\infty(x)}\,dx,\nonumber
\end{eqnarray}
$|G_t|$ is bounded by $1$ and
$p_\infty A''(F_\infty)=-\frac{2}{\sigma^2}A\times A''(F_\infty)$ is bounded,
one deduces
from \eqref{egal} that
\[
\|G_t\|^2_\infty\leq\frac{4}{\sigma^4}\|AA''\|_\infty\int_\dR\frac
{G_t^2}{p_\infty}(x)\,dx+\biggl(1\wedge\int_\dR
{\biggl(\frac{G_t}{p_\infty}(x)\biggr)'}^2p_\infty(x)\,dx\biggr).
\]
Inserting this bound in \eqref{majevol} and using Young's
inequality, one deduces that for $\eta>0$,
\begin{eqnarray*}\label{majevol2}
&&\frac{1}{2}\frac{d}{dt}\int_\dR
\frac{G_t^2}{p_\infty}(x)\,dx\\
&&\qquad\leq
\biggl(\varepsilon-\frac{\sigma^2}{2}\biggr)\int_\dR{\biggl(\frac{G_t}{p_\infty
}(x)\biggr)'}^2 p_\infty(x)\,dx\\
&&\qquad\quad{}+\frac{\|AA''\|_\infty\|A''\|^2_\infty}{4\varepsilon
\sigma^4}\biggl(\int_\dR\frac{G_t^2}{p_\infty}(x)\,dx\biggr)^2\\
&&\quad\qquad{}+\eta\biggl(1\wedge\int_\dR{\biggl(\frac{G_t}{p_\infty
}(x)\biggr)'}^2 p_\infty(x)\,dx\biggr)^2+\frac{\|A''\|^4_\infty
}{1024\varepsilon^2\eta}\biggl(\int_\dR\frac{G_t^2}{p_\infty}(x)\,dx
\biggr)^2\\
&&\qquad\leq
\biggl(\varepsilon+\eta-\frac{\sigma^2}{2}\biggr)\int_\dR{\biggl(\frac{G_t}{p_\infty
}(x)\biggr)'}^2p_\infty(x)\,dx\\
&&\quad\qquad{}+\biggl(\frac{\|AA''\|_\infty\|A''\|^2_\infty}{4\varepsilon\sigma
^4}+\frac{\|A''\|^4_\infty}{1024\varepsilon^2\eta}\biggr)\biggl(\int_\dR
\frac{G_t^2}{p_\infty}(x)\,dx\biggr)^2.
\end{eqnarray*}
One easily concludes with \eqref{poincarérep} and Lemma \ref{compedo}
below.
\end{pf}
\begin{arem}
\textup{(i)} After reading this proof, one may wonder whether one could replace
the upper bound in \eqref{majevol} by
\[
\biggl(\varepsilon-\frac{\sigma^2}{2}\biggr)\int_\dR{\biggl(\frac{G_t}{p_\infty
}(x)\biggr)'}^2 p_\infty(x)\,dx+\frac{\|A''\|^2_\infty}{16\varepsilon}\int
_\dR\frac{G_t^2}{p_\infty}(x)\,dx
\]
using $\|G_t\|_\infty\leq1$. If the constant $c$ in the Poincar\'e
inequality \eqref{poincarérep} was smaller than
$\frac{\sigma^4}{\|A''\|_\infty^2}$, one could deduce exponential
convergence of $\int_\dR\frac{G_t^2}{p_\infty}(x)\,dx$ to $0$ even for large
values of $\int_\dR\frac{G_0^2}{p_\infty}(x)\,dx$. In case
$A(x)=\frac{1}{2}(x^2-x)$ (see Example \ref{ex:explicit}), one has
$\|A''\|_\infty=1$ and
\begin{eqnarray*}
c&\geq& \int_\dR x^2p_\infty(x)\,dx-\biggl(\int_\dR
xp_\infty(x)\,dx\biggr)^2=\int_0^{+\infty}\frac{x^2}{2\sigma^2\cosh
^2({x}/{(2\sigma^2}))}\,dx\\
&>&4\sigma^4\int_0^{+\infty}y^2e^{-2y}\,dy=\sigma^4=\frac{\sigma^4}{\|A''\|
_\infty^2},
\end{eqnarray*}
and this approach does not work.

\textup{(ii)} Convexity of $A$ implies nonnegativity of the term
$A(F_t)-A(F_\infty)-A'(F_\infty)G_t$ which appears in the
right-hand side of the first displayed equality in the proof. One may
wonder if one could exploit this property
to obtain exponential convergence of $p_t$ to $p_\infty$ even if
$p_0$ is not close to $p_\infty$. We have not been able to do so.
\end{arem}
\begin{pf*}{Proof of Theorem \protect\ref{longtime}}
By \eqref{der1}, $p_\infty'=-\frac{2}{\sigma^2}A'(F_\infty)p_\infty$ and
$\|p_\infty\|_\infty\leq\frac{2\|A\|_\infty}{\sigma^2}$. The Fokker--Planck
equation \eqref{fp} for $p_t$ ensures that
\begin{eqnarray*}
&&\frac{1}{2}\frac{d}{dt}\int_\dR\frac{(p_t(x)-p_\infty(x))^2}{p_\infty
(x)}\,dx\\
&&\qquad=-\frac{\sigma^2}{2}\int_\dR
{\biggl(\frac{p_t}{p_\infty}(x)\biggr)'}^2p_\infty(x)\,dx\\
&&\qquad\quad{}-\int_\dR
\bigl(A'(F_t)-A'(F_\infty)\bigr)(x)(p_t-p_\infty)(x)\biggl(\frac{p_t}{p_\infty
}(x)\biggr)'\,dx\\
&&\qquad\quad{}-\int_\dR
\bigl(A'(F_t)-A'(F_\infty)\bigr)(x)p_\infty(x)\biggl(\frac{p_t}{p_\infty}(x)\biggr)'\,dx.
\end{eqnarray*}
Then, using Young's inequality and \eqref{majGinf}, one easily checks that
for $\varepsilon,\eta>0$,
\begin{eqnarray*}
&&\frac{1}{2}\frac{d}{dt}\int_\dR\frac{(p_t(x)-p_\infty(x))^2}{p_\infty
(x)}\,dx\\
&&\qquad\leq
\biggl(\eta+\varepsilon-\frac{\sigma^2}{2}\biggr)\int_\dR
{\biggl(\frac{p_t}{p_\infty}(x)\biggr)'}^2p_\infty(x)\,dx\\
&&\qquad\quad{}+\frac{1}{4\varepsilon}\int_\dR
\bigl(A'(F_t)(x)-A'(F_\infty)(x)\bigr)^2\frac{(p_t(x)-p_\infty(x))^2}{p_\infty
(x)}\,dx\\
&&\qquad\quad{}+\frac{1}{4\eta}\int_\dR\bigl(A'(F_t)(x)-A'(F_\infty
)(x)\bigr)^2p_\infty(x)\,dx\\
&&\qquad\leq\biggl(\eta+\varepsilon-\frac{\sigma^2}{2}\biggr)\int_\dR
{\biggl(\frac{p_t}{p_\infty}(x)\biggr)'}^2p_\infty(x)\,dx\\
&&\qquad\quad{}+\frac{\|A''\|_\infty^2}{4\varepsilon}\biggl(\int_\dR\frac
{(p_t(x)-p_\infty(x))^2}{p_\infty(x)}\,dx\biggr)^2\\
&&\qquad\quad{}+\frac{\|A''\|
_\infty^2}{4\eta}\times\frac{4\|A\|_\infty^2}{\sigma^4}\int_\dR\frac
{G_t^2}{p_\infty}(x)\,dx.
\end{eqnarray*}
By \eqref{majofoncrepdens} and Lemma \ref{convfoncrep}, for $\int_\dR
\frac{(p_0(x)-p_\infty(x))^2}{p_\infty(x)}\,dx$ small enough, the last term
of the r.h.s. is smaller than $\frac{\tilde{c} e^{-C t}}{C}\int_\dR
\frac{(p_0(x)-p_\infty(x))^2}{p_\infty(x)}\,dx$. Since $\int_\dR
{(\frac{p_t}{p_\infty}(x))'}^2p_\infty(x)\,dx$ is greater
than $
\frac{1}{c}\int_\dR\frac{(p_t(x)-p_\infty(x))^2}{p_\infty(x)}\,dx$, one
easily concludes by Lemma \ref{compedo} below.
\end{pf*}
\begin{alem}\label{compedo}
Assume that $u\dvtx\dR_+\rightarrow\dR_+$ satisfies
\[
\forall t\geq0\qquad \frac{du}{dt}(t)\leq\beta
u(t)\bigl(u(t)-\alpha\bigr)+\gamma e^{-\delta t}
\]
for some constants $\alpha,\beta,\delta>0$ and $\gamma\geq0$.

If $\gamma=0$ and $u(0)<\alpha$, then
\[
{\forall t\geq
0\qquad u(t)\leq\frac{\alpha u(0)e^{-\alpha\beta
t}}{\alpha+u(0)(e^{-\alpha\beta
t}-1)}}.
\]

If $u(0)< \frac{\alpha}{2}$ and $\gamma<\frac{\beta\alpha^2}{4}$,
then $u(t)$ converges to $0$ exponentially fast as ${t\rightarrow
+\infty}$.
\end{alem}

\begin{pf}
When $\gamma=0$, as long as $u(t)\in(0,\alpha)$, one has
\[
\frac{du}{dt}(t)\biggl(\frac{1}{u(t)}+\frac{1}{\alpha-u(t)}\biggr)\leq
-\alpha\beta
\]
and after integration one obtains the desired estimation. Since the
upper bound is not greater than $u(0)$ and $u(t)=0\Rightarrow\forall
s\geq t, u(s)=0$ one easily concludes.

Now when $\gamma\in(0,\frac{\beta\alpha^2}{4})$, one has $\beta
a(\alpha-a)=\gamma$ for some $a\in(0,\frac{\alpha}{2})$ and
\[
\frac{d}{dt}\biggl(u(t)\wedge\frac{\alpha}{2}-a\biggr)^+=1_{\{a<u(t)<{\alpha
}/{2}\}}\frac{du}{dt}(t)\leq
0.
\]
Hence when $u(0)<\frac{\alpha}{2}$, $\forall t\geq0$, $u(t)\leq
u(0)\vee a$ and
\[
\frac{du}{dt}(t)\leq-\beta\bigl(\alpha-u(0)\vee a\bigr) u(t)+\gamma
e^{-\delta t}.
\]
For $v(t)=e^{\beta(\alpha-u(0)\vee a)t}u(t)$ one deduces
\[
\frac{dv}{dt}(t)\leq\gamma e^{(\beta(\alpha-u(0)\vee
a)-\delta)t}
\]
and one concludes by integration of this inequality that $u(t)$ is
bounded by\break $C(1+t)e^{-[(\beta(\alpha-u(0)\vee a))\wedge\delta]t}$.
\end{pf}
%
\subsection[The particle system (4)]{The particle system \protect\eqref{systpart}}
Let us suppose that $A(1)=0$ and that the first-order moment associated
with the initial probability measure $m$ is defined and equal to $\bar
{x}$. As in the case of the granular media
equation considered by Malrieu \cite{mal1,mal2}, the direction
$(1,1,\ldots,1)$ is quite singular for the particle system. Indeed,
\[
d(X^{1,n}_t+\cdots+X^{n,n}_t)=\sigma\sum_{i=1}^n dB^i_t,
\]
which prevents the law of $(X^{1,n}_t,\ldots,X^{n,n}_t)$ from
converging as $t\rightarrow+\infty$. Following
\cite{mal1,mal2}, one introduces the hyperplane ${\mathcal
M}_n=\{y=(y_1,\ldots,y_n)\in{\mathbb R}^n\dvtx y_1+\cdots+y_n=n\bar{x}\}$
orthogonal to this singular direction and denotes by $\bar{P}$ the orthogonal
projection on ${\mathcal M}_n$ and by $P$ the orthogonal projection on
$\{y=(y_1,\ldots,y_n)\in{\mathbb R}^n\dvtx y_1+\cdots+y_n=0\}$. Since $\sum_{i=1}^n
a_n(i)=n(A(1)-A(0))=0$, the orthogonal projection
$(Y^{i,n}_t=\bar{x}+X^{i,n}_t-\frac{1}{n}\sum_{j=1}^nX^{j,n}_t)_{1\leq
i\leq
n}$ of the original particle system on ${\mathcal M}_n$ is a
diffusion on this hyperplane solving
%
\begin{equation}\label{projeds}
dY^{i,n}_t=\sigma\frac{n-1}{n}\,dB^i_t-\frac{\sigma}{n}\sum_{j\neq
i}dB^j_t-a_n\Biggl(\sum_{j=1}^n1_{\{Y^{j,n}_t\leq
Y^{i,n}_t\}}\Biggr)\,dt.
\end{equation}

Propagation of chaos for the projected system is a consequence of the
following estimate.
\begin{aprop}\label{chaosproj}
Assume that $A$ is convex, such that $A'$ is Lipschitz continuous with
constant $K$ and $A(1)=0$ and that the initial measure $m$ has a finite
second order moment. Then, $\forall i\in\{1,\ldots,n\}, \forall t\geq
0$,
\[
\E[(X^i_t-Y^{i,n}_t)^2]\leq
\frac{1}{n}\biggl[\frac{K^2t^2}{6}+\E[(X_0-\bar{x})^2]+\sigma
^2t+2\int_0^t\int_\dR A(F_s(x))\,dx\,ds\biggr],
\]
where $X^i$ is solution of \eqref{eq:nonlini}.
\end{aprop}
\begin{pf}
Denoting $X_1^n(t)=(X^1_t,\ldots,X^n_t)$,
$X^{n,n}_1(t)=(X^{1,n}_t,\ldots,X^{n,n}_t)$ and
$Y^{n,n}_1(t)=(Y^{1,n}_t,\ldots,Y^{n,n}_t)$, one has
\begin{eqnarray}
|X_1^n(t)-Y^{n,n}_1(t)|^2&=&|X_1^n(t)-\bar{P}X^{n,n}_1(t)|^2\\
&=&|X_1^n(t)-\bar{P}X^{n}_1(t)|^2+|\bar{P}X^{n}_1(t)-\bar
{P}X^{n,n}_1(t)|^2\nonumber\\
&\leq&\frac{1}{n}\Biggl(\sum_{i=1}^n
(X^i_t-\bar{x})\Biggr)^2+\sum_{i=1}^n(X^i_t-X^{i,n}_t)^2.
\label{estichaosproj}
\end{eqnarray}
Since $(X_t-\bar{x})^2\leq3((X_0-\bar{x})^2+\sigma^2
B_t^2+\|A'\|_\infty^2t^2)$, the variable $X_t$ is square
integrable. As
\begin{eqnarray*}
\forall x>0\qquad|(x-\bar{x})A(F_t(x))|&\leq&
\|A'\|_\infty\bigl(1-F_t(x)\bigr)(x+|\bar{x}|)\\
&\leq&\|A'\|_\infty\biggl(\frac{\E(X_t^2)}{x}+|\bar{x}|\bigl(1-F_t(x)\bigr)\biggr),
\end{eqnarray*}
one has $\lim_{x\rightarrow+\infty} (x-\bar{x})A(F_t(x))=0$. Similarly
$(x-\bar{x})A(F_t(x))$ also vanishes as $x\rightarrow-\infty$ and $\int
_\dR(x-\bar{x})A'(F_t(x))p_t(x)\,dx=-\int_\dR
A(F_t(x))\,dx$. Computing $(X_t-\bar{x})^2$ by It\^o's formula and taking
expectations, one deduces that
\[
\E\bigl((X_t-\bar{x})^2\bigr)=\E\bigl((X_0-\bar{x})^2\bigr)+\sigma^2t+2\int_0^t\int_\dR
A(F_s(x))\,dx\,ds.
\]
Moreover, by \eqref{evolesp}, $\E(X_t-\bar{x})=-A(1)t=0$. One
concludes by taking expectations in \eqref{estichaosproj} then using
Theorem \ref{chaos} and exchangeability of the particles.
\end{pf}

Let us now study the long time behavior of the projected particle
system.

\begin{thm}
Assume that the function $A$ is uniformly convex on $[0,1]$ with
constant $\alpha$ [see \eqref{unifconv}] and such that
$A(1)=0$. Then, the probability measure with density
\[
p^n_\infty(y)=\frac{1}{Z_n}e^{-{2}/{\sigma^2}\sum_{i=1}^na_n(i)y_{(i)}}
\]
with respect to the Lebesgue measure $dy$ on ${\mathcal M}_n$ is invariant
for the projected dynamics \eqref{projeds}. Here $y_{(1)}\leq
y_{(2)}\leq
\cdots\leq y_{(n)}$ denotes the increasing reordering of the
coordinates of
$y=(y_1,\ldots,y_n)$ and $Z_n=\int_{{\mathcal
M}_n}e^{-\frac{2}{\sigma^2}\sum_{i=1}^na_n(i)y_{(i)}}\,dy$.
Moreover, if $(Y^{1,n}_0,\ldots,Y^{n,n}_0)$ admits a symmetric density
$p^n_0(y)$ with respect to the Lebesgue measure on ${\mathcal M}_n$, then
for all $t\geq0$, $(Y^{1,n}_t,\ldots,Y^{n,n}_t)$ admits a symmetric density
$p^n_t(y)$ which is such that
%
\begin{eqnarray}\label{estlsystpart}
\forall t\geq0\qquad &&\int_{{\mathcal
M}_n}\biggl(\frac{p^n_t}{p^n_\infty}(x)-1\biggr)^2p^n_\infty(x)
\,dx\nonumber\\[-8pt]\\[-8pt]
&&\qquad\leq e^{-\lambda_nt}\int_{{\mathcal
M}_n}\biggl(\frac{p^n_0}{p^n_\infty}(x)-1\biggr)^2p^n_\infty(x)
\,dx\nonumber
\end{eqnarray}
where the sequence $(\lambda_n)_n$ is bounded from below by
$\frac{\alpha^2}{12^3\sigma^2}$.
\label{tlsystpart}
\end{thm}

In order to deduce long time properties of the nonlinear process from
long time properties of the projected system, it is not restrictive to
assume that $p^n_0$ is symmetric (see Remark \ref{reflsyst} to get some
intuition about this hypothesis). But the lack of uniformity in time of
the estimation given in Proposition \ref{chaosproj} is a real problem.
\begin{arem}
In case $n=2$, the process $Y_t=Y^{2,2}_t-Y^{1,2}_t$ solves the stochastic
differential equation
\[
dY_t=\sigma(dB^2_t-dB^1_t)-{\rm sgn}(Y_t)\bigl(a_2(2)-a_2(1)\bigr)\,dt
\]
and the density of $Y_t$ converges exponentially to
$\frac{a_2(2)-a_2(1)}{2\sigma^2}e^{(-({a_2(2)-a_2(1))}/{\sigma^2})|y|}$
when the density of $Y_0$ is close enough to this limit. As
$(Y^{1,2}_t,Y^{2,2}_t)=\overline x+\frac{1}{2}(-Z_t,Z_t)$, one easily deduces
exponential convergence of the density of $(Y^{1,2}_t,Y^{2,2}_t)$
on the straight line ${\mathcal M}_2$ to
$\frac{a_2(2)-a_2(1)}{\sqrt{2}\sigma^2}e^{-({a_2(2)}/{\sigma
^2})2y_{(2)}}e^{({a_2(1)}/{\sigma^2})(-2y_{(1)})}$.
\end{arem}

The proof of Theorem \ref{tlsystpart} relies on the following
Poincar\'e inequality.
\begin{aprop}\label{poincan}
Under the assumptions of Theorem \textup{\ref{tlsystpart}}, the density
\[
\tilde{p}^n_\infty(y)=\frac{n!1_{\{y_1\leq y_2\leq\cdots\leq
y_n\}}}{Z_n}e^{-({2}/{\sigma^2})\sum_{i=1}^na_n(i)y_i}
\]
on ${\mathcal M}_n$ is such that for $f\dvtx\dR^n\rightarrow\dR$
regular enough,
%
\begin{eqnarray}\label{poinca}
&&\int_{{\mathcal M}_n}\biggl(f(y)-\int_{{\mathcal M}_n}f(y)\tilde
{p}^n_\infty(y)\,dy\biggr)^2\tilde{p}^n_\infty(y)\,dy\nonumber\\[-8pt]\\[-8pt]
&&\qquad\leq
\frac{\sigma^2}{\lambda_n}\int_{{\mathcal M}_n}|P\nabla f(y)|^2\tilde
{p}^n_\infty(y)\,dy\nonumber
\end{eqnarray}
where the sequence $(\lambda_n)_n$ is bounded from below by
$\frac{\alpha^2}{12^3\sigma^2}$.
\end{aprop}
\begin{pf*}{Proof of Theorem \protect\ref{tlsystpart}}
Let us first check the following Green formula: for
$f\dvtx\dR^n\rightarrow\dR$ and $u\dvtx\dR^n\rightarrow\dR^n$ regular
enough,
%
\begin{equation}\label{green}
\int_{{\mathcal M}_n}f\nabla\cdot(Pu)(y)\,dy=-\int_{{\mathcal
M}_n}P\nabla
f\cdot(Pu)(y)\,dy.
\end{equation}
Let ${\mathbf1}\in\dR^n$ denote the vector with all coordinates equal
to $1$.
For $\varphi\dvtx\dR\rightarrow\dR$ and $v\dvtx\dR^n\rightarrow\dR^n$, one has
\begin{eqnarray*}
&&\int_{\dR}\varphi\bigl(\sqrt{n}z\bigr)\int_{{\mathcal M}_n}\nabla\cdot(Pv)
\biggl(y+\frac{z{\mathbf1}}{\sqrt{n}}\biggr)\,dy\,dz\\
&&\qquad=\int_{\dR^n}\varphi(x_1+\cdots+x_n-n\bar{x})\nabla\cdot(Pv)(x)\,dx\\
&&\qquad=-\int_{\dR^n}\varphi'(x_1+\cdots+x_n-n\bar{x}){\mathbf1}\cdot
(Pv)(x)\,dx=0.
\end{eqnarray*}
The function $\varphi$ being arbitrary, one deduces that
$\int_{{\mathcal M}_n}\nabla\cdot(Pv)(y)\,dy=0$.\vspace*{2pt} Since $\nabla\cdot
P(fu)=\nabla
f\cdot(Pu)+f\nabla\cdot(Pu)=P\nabla
f\cdot(Pu)+f\nabla\cdot(Pu)$, \eqref{green} follows for the choice
$v=fu$.

By weak uniqueness for \eqref{projeds}, when
$(Y^{1,n}_0,\ldots,Y^{n,n}_0)$ has a symmetric density $p^n_0$ with
respect to the Lebesgue measure on ${\mathcal M}_n$, the particles
$Y^{i,n}$, $i\in\{1,\ldots,n\}$ are exchangeable and for each $t\geq
0$, $(Y^{1,n}_t,\ldots,Y^{n,n}_t)$ has a symmetric density $p^n_t$. By
composition with the projection $\bar{P}$, one obtains an extension of
$p^n_t$ on $\dR^n$ that we still denote by $p^n_t$. Since
$\sum_{i=1}^n a_n(i)=n(A(1)-A(0))=0$, setting
\[
b(y)=\sum_{\tau\in{\mathcal S}_n}1_{\{y_{\tau(1)}\leq
y_{\tau(2)}\leq\cdots\leq
y_{\tau(n)}\}}
\pmatrix{{c}a_n(\tau^{-1}(1)) \cr
a_n(\tau^{-1}(2)) \cr
\vdots\cr
a_n(\tau^{-1}(n))
},
\]
one has $Pb=b$ and the infinitesimal generator associated with \eqref{projeds}
is $L\psi=\frac{\sigma^2}{2}\nabla\cdot(P\nabla\psi)-Pb\cdot\nabla\psi$.
Computing $d\psi(Y^{1,n}_t,\ldots,Y^{n,n}_t)$ by It\^o's formula and taking
expectations then using \eqref{green}, one obtains
\begin{eqnarray*}
\int_{{\mathcal M}_n}\psi(y)\partial_tp^n_t(y)\,dy&=&\int
_{{\mathcal
M}_n}L\psi(y) p^n_t(y) \,dy\\
&=&\int_{{\mathcal
M}_n}\psi(y)\nabla\cdot P\biggl(\frac{\sigma^2}{2}\nabla
p^n_t+bp^n_t\biggr)(y)\,dy.
\end{eqnarray*}
Hence the densities solve the Fokker--Planck equation
\[
\partial_t p^n_t=\nabla\cdot P\biggl(\frac{\sigma^2}{2}\nabla
p^n_t+bp^n_t\biggr).
\]
Now using \eqref{green} and $b=-\frac{\sigma^2\nabla
p^n_\infty}{2p^n_\infty}$, one deduces
\begin{eqnarray}\label{evolpnt}
&&\partial_t\int_{{\mathcal
M}_n}\biggl(\frac{p^n_t}{p^n_\infty}(y)-1\biggr)^2p^n_\infty(y) \,dy\nonumber\\
&&\qquad=2\int_{{\mathcal M}_n}\frac{p^n_t}{p^n_\infty}(y)\nabla\cdot
P\biggl(\frac{\sigma^2}{2}\nabla p^n_t+bp^n_t\biggr)(y)\,dy\nonumber\\[-8pt]\\[-8pt]
&&\qquad=-\sigma^2\int_{{\mathcal
M}_n}P\nabla\frac{p^n_t}{p^n_\infty}(y)\cdot P\frac
{\nabla
p^n_t+({2bp^n_t}/{\sigma^2})}{p^n_\infty}(y)p^n_\infty(y) \,dy\nonumber\\
&&\qquad=-\sigma^2\int_{{\mathcal
M}_n}\bigg|P\nabla\frac{p^n_t}{p^n_\infty}(y)\bigg|^2p^n_\infty(y)\, dy.\nonumber
\end{eqnarray}
By symmetry of the function $\frac{p^n_t}{p^n_\infty}$ and \eqref{poinca},
\begin{eqnarray*}
\sigma^2\int_{{\mathcal
M}_n}\bigg|P\nabla\frac{p^n_t}{p^n_\infty}(y)\bigg|^2p^n_\infty(y)
\,dy&=&\sigma^2\int_{{\mathcal
M}_n}\bigg|P\nabla\frac{p^n_t}{p^n_\infty}(y)\bigg|^2\tilde
{p}^n_\infty(y) \,dy\\
&\geq&\lambda_n \int_{{\mathcal
M}_n}\biggl(\frac{p^n_t}{p^n_\infty}(y)-1\biggr)^2\tilde
{p}^n_\infty(y) \,dy\\
&\geq&\lambda_n\int_{{\mathcal
M}_n}\biggl(\frac{p^n_t}{p^n_\infty}(y)-1\biggr)^2p^n_\infty(y) \,dy
\end{eqnarray*}
and the conclusion follows.
\end{pf*}

Notice that the computation in \eqref{evolpnt} is formal and can only be
justified when $p^n_t$ is a smooth solution of the Fokker--Planck equation.
\begin{arem}\label{reflsyst}
Let us denote by $Y^{(1),n}_t\leq\cdots\leq Y^{(n),n}_t$ the
increasing reordering of $(Y^{1,n}_t,\ldots,Y^{n,n}_t)$. According
to \cite{jou2002}, the reordered system is a diffusion process
normally reflected at the boundary of the closed convex set $\{y\in
{\mathcal M}_n\dvtx y_1\leq y_2\leq\cdots\leq y_n\}$. More precisely,
%
\begin{equation}
\cases{
dY^{(i),n}_t=\sigma
\,d\beta^i_t-a_n(i)\,dt+(\gamma^i_t-\gamma^{i+1}_t)\,d|K|_t,\cr
\biggl(\displaystyle\int_0^t(\gamma^i_s-\gamma^{i+1}_s)\,d|K|_s,1\leq i\leq n\biggr)_{t\geq
0}\mbox{ is a continuous process}\vspace*{2pt}\cr
\qquad\mbox{with finite variation equal to
}|K|_t, \cr
\gamma^1\equiv\gamma^{n+1}\equiv0,\cr
d|K|_t\mbox{ \textup{a.e.} }\forall2\leq i\leq n,  \gamma^i_t\geq
0\mbox{ and }\gamma^i_t\bigl(Y^{(i),n}_t-Y^{(i-1),n}_t\bigr)=0,
}
\end{equation}
where $(\beta^1,\ldots,\beta^n)$ is a Brownian motion such that
$\frac{{\langle\beta^i,\beta^j\rangle}_t}{t}=1_{\{i=j\}}-1/n$.

If the initial condition $(Y^{(1),n}_0\leq\cdots\leq Y^{(n),n}_0)$
admits a density $\tilde{p}^n_0$ with respect to the Lebesgue measure
on ${\mathcal M}_n$, then the law of
$(Y^{(1),n}_t,\ldots,Y^{(n),n}_t)$ is the image by increasing reordering
of the symmetric law
of the solution $(Y^{1,n}_t,\ldots,Y^{n,n}_t)$ to \eqref{projeds}
starting from $(Y^{1,n}_0,\ldots,Y^{n,n}_0)$ with density $p^n_0$
obtained by symmetrization of $\tilde{p}^n_0$. Therefore
$(Y^{(1),n}_t,\ldots,Y^{(n),n}_t)$ has the density
$\tilde{p}^n_t(y)=n!p^n_t(y)1_{\{y_1\leq\cdots\leq y_n\}}$ and
\eqref{estlsystpart} holds with $p^n$ replaced by
$\tilde{p}^n$.
\end{arem}

In order to prove Proposition \ref{poincan},
we take advantage of the specific form of the density
$\tilde{p}^n_\infty$. Remarking that $\tilde{p}^n_\infty$ is the
density of the image of a vector of independent exponential random
variables by a linear transformation, one first obtains the following
result.

\begin{alem}\label{poincal}
The Poincar\'e inequality \eqref{poinca} holds with the constant
$\lambda_n$ greater than $\frac{\alpha^2}{4\sigma^2}$ multiplied by
the  smallest\vspace*{2pt} eigenvalue $\tilde{\lambda}_n$ of the
$(n-1)\times(n-1)$ matrix $Q^n$ defined by $\forall1\leq i,j\leq
n-1$, $Q^n_{ij}=b_n(i)L^n_{ij}b_n(j)$ where
\[
b_n(i)=\frac{i(n-i)}{n}\quad\mbox{and}\quad L^n=\pmatrix{
2 & -1 & 0 &\ldots& \ldots& \ldots& 0\cr
-1 & 2 & -1 & 0 & \ldots& \ldots& 0\cr
0&-1&2&-1& 0 & \ldots& 0\cr
\vdots&\vdots&\vdots&\vdots&\vdots&\vdots&\vdots\cr 0&\ldots
&0&-1&2&-1&0\cr0 &\ldots& \ldots& 0 & -1&2&-1\cr
0 &\ldots& \ldots& \ldots& 0&-1&2
}
.
\]
\end{alem}

The last statement in Proposition \ref{poincan} then follows from the
next lemma which is obtained by interpreting $Q^n$ as a finite element
rigidity matrix associated with the operator
$-x(1-x)\partial_{xx}(x(1-x).)$ acting on functions on $(0,1)$. The
Hardy inequality stated in Lemma \ref{lemhard} ensures that it is
enough to bound the smallest eigenvalue of the corresponding mass
matrix from below. The resort to this one-dimensional Poincar\'e-like
inequality in order to estimate the constant in the \mbox{$n$-dimensional}
Poincar\'e inequality \eqref{poinca} is striking.

\begin{alem}\label{minolam}
The sequence $(\tilde{\lambda}_n)_n$ is bounded from below by
$1/(16\times27)$.
\end{alem}
\begin{pf*}{Proof of Lemma \protect\ref{poincal}}
Let $f$ be such that $\int_{{\mathcal
M}_n}f(y)\tilde{p}^n_\infty(y)\,dy=0$. Since the left-hand side in
the Poincar\'e inequality \eqref{poinca} only depends on the
restriction of $f$ to ${\mathcal M}_n$, one may assume that $\forall
x\in\dR^n, f(x)=f(\bar{P}x)$, which ensures that for
$(x_1,\ldots,x_n)\in\dR^n$ such that $x_1+\cdots+x_n=0$,
$f(\bar{x}+x_1,\ldots,\bar{x}+x_n)=f(x_1,\ldots,x_n)$ and $P\nabla
f(\bar{x}+x_1,\ldots,\bar{x}+x_n)=\nabla f(x_1,\ldots,x_n)$. Therefore
the Poincar\'e inequality \eqref{poinca} is
equivalent to $I(f)\leq\frac{\sigma^2}{\lambda_n} I(|\nabla f|)$
where
\[
I(g)=\int_{\dR^{n-1}}(g^2\tilde{p}^n_\infty)\bigl(-(x_2+\cdots
+x_n),x_2^n\bigr)\,dx_2^n \qquad\mbox{with } x_2^n=(x_2,\ldots,x_n).
\]
To integrate the coordinates over
independent domains, we make the change of variables $z_2^n=Mx_2^n$ where
\[
M=\pmatrix{
2 & 1 & 1 &\ldots& \ldots& 1\cr
-1 & 1 & 0 &\ldots& \ldots& 0\cr
0&-1&1 &0& \ldots& 0\cr
\vdots&\vdots&\vdots&\vdots&\vdots&\vdots\cr
0 &\ldots& 0 & -1 & 1 &0\cr
0 &\ldots& \ldots& 0&-1&1
}
.
\]
One easily checks that for $2\leq i\leq n$,
$z_2+\cdots+z_i=x_2+\cdots+x_n+x_i$ and deduce that
$(n-1)z_2+(n-2)z_3+\cdots+2z_{n-1}+z_{n}=n(x_2+\cdots+x_n)$.
Therefore
\[
M^{-1}=\frac{1}{n}\pmatrix{ 1 & 2-n & 3-n &
4-n & \ldots
& -1 \cr
1 & 2 & 3-n&4-n& \ldots& -1\cr
1&2&3 & 4-n& \ldots&-1\cr
\vdots&\vdots&\vdots&\vdots&\vdots&\vdots\cr1 &2& 3 & \ldots
&n-2&-1\cr
1 &2&3 & \ldots&\ldots&n-1
}
\]
and denoting
\[
N=\pmatrix{
 \dfrac{1-n}{n} & \dfrac{2-n}{n} &
\ldots
-\dfrac{2}{n} & -\dfrac{1}{n} \cr
& & M^{-1}& &
}
,
\]
one has
\[
I(f)=\frac{n!}{Z_n}\int_{(\dR_+)^{n-1}}f^2(Nz_2^n)e^{(-{2}/{\sigma
^2})\sum_{i=2}^{n}\beta_n(i)z_{i}}\frac{dz_2^n}{|M|}
\]
where
\begin{eqnarray*}
\beta_n(i)&=&\frac{1}{n}\bigl[(i-1)\bigl(a_n(i)+\cdots+a_n(n)\bigr)\\
&&\phantom{\frac{1}{n}[}{}-(n+1-i)\bigl(a_n(1)+\cdots+a_n(i-1)\bigr)\bigr]\\
&=&-nA\bigl((i-1)/n\bigr)>0.
\end{eqnarray*}
Here $|M|$
denotes the determinant of the matrix $M$; it is equal to $n$ by an easy
computation. The one-dimensional exponential density with parameter $c$
satisfies the Poincar\'e inequality with optimal constant $4/c^2$.
Tensorizing this inequality (see Chapters 3 and 6 in
\cite{abcfgmrs} for further details), one obtains
\begin{eqnarray*}
I(f)&\leq&\frac{n!}{Z_n}
\int_{(\dR_+)^{n-1}}\sum_{j=2}^n\frac{\sigma^4}{\beta^2_n(j)}\Biggl(\sum
_{k=1}^nN_{kj-1}\partial_kf(Nz_2
^n)\Biggr)^2e^{(-{2}/{\sigma^2})\sum_{i=2}^{n}\beta_n(i)z_{i}}\frac
{dz_2^n}{|M|}\\
&=&\int_{\dR^{n-1}}\sum_{k,l=1}^n\sum_{j=2}^n\frac{\sigma^4}{\beta
^2_n(j)}N_{kj-1}N_{lj-1}\partial_kf\partial_lf\tilde{p}^n_\infty
\bigl(-(x_2+\cdots+x_n),x_2^n\bigr)\,dx_2^n.
\end{eqnarray*}
Since $A$ is uniformly convex with constant $\alpha$ and $A(0)=A(1)=0$,

\[
\beta_n(i)=-nA\bigl((i-1)/n\bigr)\geq-\frac{n\alpha}{2}\times\frac
{i-1}{n}\biggl(\frac{i-1}{n}-1\biggr)=\frac{\alpha
}{2}b_n(i-1).
\]
Therefore
\begin{eqnarray*}
I(f)&\leq&
\frac{4\sigma^4}{\alpha^2}\int_{\dR^{n-1}}\sum_{k,l=1}^n\sum
_{j=1}^{n-1}\frac{N_{kj}N_{lj}}{b^2_n(j)}\partial_kf\partial_lf\tilde
{p}^n_\infty\bigl(-(x_2+\cdots+x_n),x_2^n\bigr)\,dx_2^n\\
&\leq&\frac{4\sigma^2}{\alpha^2\tilde{\lambda}_n}I(|\nabla
f|)
\end{eqnarray*}
where $\tilde{\lambda}_n$ denotes the inverse of the largest eigenvalue
of the
symmetric positive semidefinite matrix $\bar{N}\bar{N}^*$ defined by
$\bar{N}_{ij}=\frac{N_{ij}}{b_n(j)}$. To prove Proposition
\ref{poincan} with a possibly modified lower bound, it is enough to
check that the largest eigenvalue is bounded from above uniformly in
$n$. Unfortunately, the trace of the matrix can be bounded from below
by a positive constant multiplied by $\log(n)$. Therefore one has to
be more precise.

Let $w$ be an eigenvector associated with the largest eigenvalue:
$\bar{N}\bar{N}^*w=\frac{1}{\tilde{\lambda}_n} w$. Of course
$\bar{N}^*w$ is nonzero and multiplying the previous equality by
$\bar{N}^*$, one obtains that $\bar{N}^*w$ is an eigenvector of
$\bar{N}^*\bar{N}$ associated with the eigenvalue
$\frac{1}{\tilde{\lambda}_n}$. By symmetry,
$\frac{1}{\tilde{\lambda}_n}$ is also the largest eigenvalue of
$\bar{N}^*\bar{N}$. We are going to check that the latter matrix is
invertible with inverse equal to $Q^n$ in order to conclude the proof.
Because of the definition of $\bar{N}$, it is enough to check that
$N^*N$ is invertible with inverse equal to
$L_n$.

By construction of the matrix $N$, for the equation $Nz_2^n=x$ where
$x\in\dR^n$ to have a solution $z_2^n$, it is necessary and sufficient
that
$x_1=-(x_2+\cdots+x_n)$ and then $z_2^n=Mx_2^n$.

Now for fixed $y\in\dR^{n-1}$, let us find $x_2^n\in\dR^{n-1}$ such
that $N^*x=y$ where $x=-(x_2+\cdots+x_n,x_2^n)$. This equation writes
\[
\left((M^{-1})^*-
\pmatrix{
N_{11} &N_{11} &\ldots&N_{11}\cr
N_{12}&N_{12}&\ldots&N_{12}\cr\vdots&\vdots&\vdots&\vdots\cr
N_{1n-1}&N_{1n-1}&\ldots&N_{1n-1}
}
\right)x_2^n=y.
\]
One easily checks that the $(n-1)\times(n-1)$ matrix in the
left-hand side is equal to
\[
\pmatrix{
1&1&1&\ldots&1\cr0&1&1&\ldots&1\cr
0&0&1&\ldots&1\cr\vdots&\vdots&\vdots&\vdots&\vdots\cr0&\ldots&0&1&1\cr
0&\ldots&0&0&1
}
\qquad\mbox{with inverse
}R=\pmatrix{
1&-1&0&0&\ldots&0\cr0&1&-1&0&\ldots&0\cr\vdots&\vdots
&\vdots&\vdots&\vdots&\vdots\cr0&\ldots&0&1&-1&0\cr0&\ldots&0&0&1&-1\cr
0&\ldots&0&0&0&1
}
.
\]
Combining $x_2^n=Ry$ with the solution of the previous problem,
one obtains that the unique solution of the equation $N^*Nz_2^n=y$ is
$z_2^n=MRy$. One concludes by checking that the matrix $MR$ is equal
to $L_n$.
\end{pf*}
\begin{pf*}{Proof of Lemma \protect\ref{minolam}}
For $i\in\{1,\ldots,n-1\}$, the functions
\[
u_i(x)=
\cases{
0,&\quad if $x\in(0,1)\Big\backslash\biggl[\dfrac{i-1}{n},\dfrac{i+1}{n}\biggr]$,\cr
\dfrac{i(n-i)(x-{(i-1)}/{n})}{\sqrt{n} x(1-x)},&\quad if $x\in\biggl[\dfrac
{i-1}{n},\dfrac{i}{n}\biggr]$,\cr
\dfrac{i(n-i)({(i+1)}/{n}-x)}{\sqrt{n} x(1-x)},&\quad if $x\in\biggl[\dfrac
{i}{n},\dfrac{i+1}{n}\biggr]$,
}
\]
are such that
\[
\forall
i,j\in\{1,\ldots,n-1\}\qquad Q^n_{ij}=\int_0^1\bigl(x(1-x)u_i(x)\bigr)'\bigl(x(1-x)u_j(x)\bigr)'\,dx.
\]
By the Hardy inequality stated in Lemma \ref{lemhard} below, the
smallest eigenvalue of the matrix $Q^n$ is greater than the smallest
eigenvalue of the $(n-1)\times(n-1)$
tridiagonal matrix $R^n_{ij}=\int_0^1u_i(x)u_j(x)\,dx$ divided by $16$.

For $i\in\{1,\ldots,n-2\}$, let
$r^n_i=\int_{i/n}^{(i+1)/n}u_i(u_i-u_{i+1})(x)\,dx$ and
\[
r^n_{n-1}=\int_{(n-1)/n}^1u^2_{n-1}(x)\,dx=\frac{(n-1)^2}{n}\int
_{(n-1)/n}^1\frac{1}{x^2}\,dx=\frac{n-1}{n}.
\]
Using the change of variables $y=1-x$, one easily checks that
\[
\forall
i\in\{1,\ldots,n-1\}\qquad
R^n_{ii}-R^n_{ii-1}-R^n_{ii+1}=r^n_i+r^n_{n-i},
\]
where by convention $R^n_{10}=R^n_{n-1n}=0$. We are going to prove
that
\[
\forall n\geq3\ \forall
i\in\{2,\ldots,n-3\}\qquad r^n_i\geq\tfrac{1}{27},
\]
and that $r^n_1$ and $r^n_{n-2}$ are nonnegative. For $y\in\dR^{n-1}$,
one deduces that
\begin{eqnarray*}
y^*R^ny&=&\sum_{i=1}^{n-1}R^n_{ii}y_i^2+2\sum
_{i=1}^{n-2}R^n_{ii+1}y_iy_{i+1}\\
&=&\sum_{i=1}^{n-1}(R^n_{ii}-R^n_{ii-1}-R^n_{ii+1})y_i^2+\sum
_{i=1}^{n-2}R^n_{ii+1}(y_i+y_{i+1})^2\geq\frac{|y|^2}{27}
\end{eqnarray*}
and the conclusion follows.

Let us first suppose that $i\leq\lfloor\frac{n}{2}\rfloor-1$, which
ensures that the function $f(x)=x^2(1-x)^2$ is increasing on
$[i/n,(i+1)/n]$. Let $g(x)=u_i(u_i-u_{i+1})(x)$. One easily checks
that
\begin{eqnarray*}
\int_{i/n}^{(i+1)/n}g(x)\,dx&=&\frac{i^2(n-i)^2}{n^4}\biggl(\frac
{1}{3}-\frac{(i+1)(n-i-1)}{6i(n-i)}\biggr)\\
&\geq&
\cases{
0, & \quad if $i=1$,\cr
\dfrac{i^2(n-i)^2}{12n^4},&\quad if $i\geq2$.
}
\end{eqnarray*}
Since there is some $x_i\in[i/n,(i+1)/n]$ such that the function
$g(x)$ is non-negative on $[i/n,x_i]$ then nonpositive on
$[x_i,(i+1)/n]$, and $f$ is positive and increasing, one deduces that
for all $x\in[i/n,(i+1)/n]$, $\int_{i/n}^x\frac{g(y)}{f(y)}\,dy\geq0$.
This ensures that $\forall x\in[i/n,(i+1)/n]$
\[
\frac{d}{\,dx}\biggl(f(x)\int_{i/n}^x\frac{g(y)}{f(y)}\,dy\biggr)
=f'(x)\int_{i/n}^x\frac{g(y)}{f(y)}\,dy+g(x)\geq g(x).
\]
Therefore
\begin{eqnarray*}
r^n_i&=&\int_{i/n}^{(i+1)/n}\frac{g(y)}{f(y)}\,dy\geq
\frac{1}{f((i+1)/n)}\int_{i/n}^{(i+1)/n}g(y)\,dy\\
&\geq&
\cases{
0,& \quad if $i=1$,\cr
\dfrac{i^2(n-i)^2}{12(i+1)^2(n-i-1)^2}\geq\dfrac{1}{27},&\quad if $i\geq2$.
}
\end{eqnarray*}
Let us now suppose that $i\geq\lfloor\frac{n+1}{2}\rfloor$ so that
the function $f$ is decreasing on $[i/n,(i+1)/n]$. We deduce that
\begin{eqnarray*}
r^n_i&\geq&
\frac{1}{f(i/n)}\int_{i/n}^{(i+1)/n}(fu_i^2)(x)\,dx\\
&&{}-\frac
{1}{f((i+1)/n)}\int_{i/n}^{(i+1)/n}(fu_iu_{i+1})(x)\,dx\\
&=&\frac{1}{3}-\frac{i(n-i)}{6(i+1)(n-i-1)}
\end{eqnarray*}
and the left-hand side is greater than $1/12$ for $i\leq n-3$ and
nonnegative for $i=n-2$.

We still have to deal with the case $n$ odd and $i=(n-1)/2$. Then, $f$
is not monotonic on $I_n=[i/n,(i+1)/n]=[1/2-1/2n,1/2+1/2n]$. But by
symmetry,
\begin{eqnarray*}
r^n_{(n-1)/2}&=&\frac{(n-1)^2(n+1)^2}{16n}\int
_{1/2-1/2n}^{1/2+1/2n}\frac{(1/2+1/2n-x)(1-2x)}{x^2(1-x)^2}\,dx\\
&=&\frac
{(n-1)^2(n+1)^2}{32n}\int_{1/2-1/2n}^{1/2+1/2n}\frac
{(1-2x)^2}{x^2(1-x)^2}\,dx\\
&\geq&\frac{(n-1)^2(n+1)^2}{2n}\int
_{1/2-1/2n}^{1/2+1/2n}(1-2x)^2\,dx=\frac{(n^2-1)^2}{6n^4},
\end{eqnarray*}
which completes the proof.
\end{pf*}

\begin{alem}\label{lemhard}
For all $u\in L^2(0,1)$ such that the distribution derivative
${(x(1-x)u(x))'}$ belongs to $L^2(0,1)$,
\[
\int_0^1u^2(x)\,dx\leq16\int_0^1\bigl(\bigl(x(1-x)u(x)\bigr)'\bigr)^2\,dx.
\]
\end{alem}
\begin{pf}
For $v$ a $C^\infty$ function with compact support on $(0,1)$, by
the integration by parts formula,
\begin{eqnarray*}
\int_0^{1/2}\frac{v^2(x)}{x^2(1-x)^2}\,dx&\leq&
4\displaystyle\int_0^{1/2}\frac{v^2(x)}{x^2}\,dx=
8\biggl(\int_0^{1/2}\frac{vv'(x)}{x}\,dx-v^2(1/2)\biggr)\\&\leq&8
\biggl(\displaystyle\int_0^{1/2}\frac{v^2(x)}{x^2}\,dx\biggr)^{1/2}\biggl(\int
_0^{1/2}(v'(x))^2\,dx\biggr)^{1/2}.
\end{eqnarray*}
Dealing with the integral on $(1/2,1)$ in a symmetric way, one deduces
%
\begin{equation}
\int_0^{1}\frac{v^2(x)}{x^2(1-x)^2}\,dx\leq16\int
_0^{1}(v'(x))^2\,dx.\label{hardy}
\end{equation}
Now approximating $v\in H^1_0(0,1)$ by a sequence of $C^\infty$
functions with compact support converging in the $H^1$ norm and almost
everywhere, one deduces with the Fatou lemma that the inequality still
holds
for $v\in H^1_0$.

For $u$ satisfying the hypotheses in the lemma, $v(x)=x(1-x)u(x)$
belongs to $H^1(0,1)$. According to Theorem VIII.2, page 122 of
\cite{brezis}, $v$ admits a representative continuous on $[0,1]$ still
denoted by $v$. Moreover, since $u(x)=\frac{v(x)}{x(1-x)}$ belongs to
$L^2(0,1)$, necessarily, $v(0)=v(1)=0$. By Theorem VIII.11, page 133
of \cite{brezis}, $v$~belongs to
$H^1_0(0,1)$ and the conclusion follows from \eqref{hardy}.
\end{pf}
\section*{Acknowledgment}
We warmly thank Tony Leli\`evre (CERMICS) for fruitful
discussions concerning the analysis of the long time behavior of the
nonlinear process.

%

%
\printaddresses
\end{document}